\newcommand{\ddatt}{                                   %%%
}%%                                                    %%%
\documentclass[12pt]{article}
\usepackage{euscript,amsmath, amssymb, amsfonts}
\usepackage{color}

\newcommand{\ig}{\Big}
\newcommand{\igg}{\Bigg}

\newcommand{\sigmama}{ {{g}} }

\newcommand{\p}{\partial}
\newcommand{\R}{{\cal R}}
\newcommand{\HH}{ H_{W(\cal R)}(\nu) }
\newcommand{\G}{{\open C}[W(\R)]}

\newcommand{\str}{sp }
\newcommand{\be}{\begin{equation}}
\newcommand{\ee}{\end{equation}}
\newcommand{\bee}{\begin{eqnarray}}
\newcommand{\eee}{\end{eqnarray}}
\newcommand\nn{\nonumber \\}
\newcommand\n{ }

\newcommand{\x}{{\vec x}}
\newcommand{\yy}{{\vec y}}
\newcommand{\cc}{{\vec c}}
\newcommand{\vv}{{\vec v}}

\newcounter{theorem}

\newenvironment{proof}[1][Proof]{\noindent\textsf{#1.\ }}{ \ \rule{0.5em}{0.5em}}

\renewcommand{\theequation}{\arabic{equation}}
\font\frtnfr=eufm10   scaled\magstep1
\font\twlfr=eufm10
\font\tenfr=eufm10

\newfam\frfam
\textfont\frfam=\frtnfr
\scriptfont\frfam=\twlfr
\scriptscriptfont\frfam=\tenfr
\def\fr{\fam\frfam}

\font\frtnopen=msbm10  scaled\magstep2
\font\twlopen=msbm10
\font\tenopen=msbm10

\newfam\openfam
\textfont\openfam=\frtnopen
\scriptfont\openfam=\twlopen
\scriptscriptfont\openfam=\tenopen
\def\open{\fam\openfam}

\font\frtnsf = cmss12 scaled\magstep1
\font\twlsf = cmss10
\font\tensf = cmss9

\newfam\Scfam
\textfont\Scfam = \frtnsf
\scriptfont\Scfam = \twlsf
\scriptscriptfont\Scfam = \tensf

\begin{document}

%%%%%%%%%%%%%%%%%%%%%%%%%%%%%%%%%%%%%%%%%%%%%%%%%%%%%%%%%

%%%%%%%%%%%%%%%%%%%%%%%%%%%%%%%%%%%%%%%%%%%%%%%%%%%%%%%%%%%%%%%%%%%%%
%%%%%%%%%%%%%%      TITLE PAGE     %%%%%%%%%%%%%%%%%%%%%%%%%%%%%%%%%%
%%%%%%%%%%%%%%%%%%%%%%%%%%%%%%%%%%%%%%%%%%%%%%%%%%%%%%%%%%%%%%%%%%%%%

\sloppy
\title
 {
                  \vspace{1cm}
      Traces on the Algebra of Observables of
      Rational Calogero Model based on the Root System
 }
\author
{ S.E.Konstein\thanks{E-mail: konstein@lpi.ru}\ \ and
 I.V.Tyutin\thanks{E-mail: tyutin@lpi.ru
}          \thanks
             {This work is supported by the Russian Fund for Basic
               Research,
              Grant 11-02-00685.}
  \\
           {\it {\small} I.E.Tamm Department of Theoretical Physics,}
  \\
               {\it {\small} P. N. Lebedev Physical Institute,}
  \\
         {\it {\small} 117924, Leninsky Prospect 53, Moscow, Russia.}
 }
\date{
\ddatt 
}
%--------------------------------------------------------------------
\maketitle
%--------------------------------------------------------------------
%%%%%%%%%%%%%%%%%%%%%%%%%%%%%%%%%%%%%%%%%%%%%%%%%%%%%%

\begin{abstract}
It is shown that$H_{W(\R)}(\nu)$, the algebra  of observables of the
rational Calogero model based on the root system $\R \subset {\open R}^N$, possesses
$T_\R$ independent traces, where $T_\R$ is the number of conjugacy
classes of elements without eigenvalue $1$ belonging to the Coxeter
group $W(\R)\subset End({\open R}^N)$ generated by the root system $\R$.

Simultaneously, we reproduced an older result:
the algebra $H_{W(\R)}(\nu)$,  considered as a superalgebra with a natural parity,
possesses
$ST_\R$ independent  supertraces, where $ST_\R$ is the number of conjugacy
classes of elements without eigenvalue $-1$ belonging
to $W(\R)$.
\end{abstract}
%%%%%%%%%%%%%%%%%%%%%%%%%%%%%%%%%%%%%%%%%%%%%%%%%%%%%%

\section{Introduction}

It was shown in \cite{KV} and \cite{Ko}
that for every associative superalgebra $H_{W(\R)}(\nu)$ of observables
of the rational Calogero model based on the root system $\R$, the space of
supertraces is nonzero. The dimensions of these spaces for every root system are listed in
\cite{stek}.

Here we also consider  these superalgebras as algebras (parity forgotten%
)
and find the
conditions of existence and the dimensions of the spaces of traces on these algebras.

Astonishingly, the proof differs from the one in \cite{KV} and \cite{Ko} in several signs only,
and we provide it here indicating change of signs by means of a parameter  $\varkappa$ with
$\varkappa=-1$ for the supertraces and $\varkappa=+1$ for the traces. 
As a result, some parts of this 
text are almost copypasted from \cite{KV} and \cite{Ko}, especially Subsection \ref{genrel} and
Appendices.

\subsection{Main definitions}

\subsubsection{Traces}

Let ${\cal A}$ be an associative superalgebra with parity $\pi$.
All expressions of linear algebra are given for homogenious elements only
and are supposed to be extended to inhomogeneous elements via linearity.

%\begin{definition}\label{str}
A linear function $str$ on ${\cal A}$ is called a {\it supertrace} if
$$str(fg)=(-1)^{\pi(f)\pi(g)}str(gf) \ \mbox{ for all } f,g\in {\cal A}.$$
%\end{definition}

%\begin{definition}\label{tr}
A linear function $tr$ on ${\cal A}$ is called a {\it trace} if
$$tr(fg)=tr(gf) \ \mbox{ for all } f,g\in {\cal A}.$$
%\end{definition}

Let $\varkappa=\pm 1$.
We can unify the definitions of trace and supertrace by introducing a $\varkappa$-trace.

%\begin{definition}
We say that a linear function $\str$
\footnote{From German word {\it spur}.}
on ${\cal A}$ is a {\it $\varkappa$-trace} if
\begin{equation}\label{scom}
\str(fg)=\varkappa^{\pi(f)\pi(g)}\str(gf) \ \mbox{ for all $f,g\in {\cal A}$}.
\end{equation}
%\end{definition}

%\begin{definition}
 A linear function $L$ is {\it even} if $L(f)=0$
for any $f\in{\cal A}$ such that $\pi(f)=1$,
and is
{\it odd} if $L(f)=0$
for any $f\in{\cal A}$ such that $\pi(f)=0$.
%\end{definition}

\medskip
Let ${\cal A}_1$ and
${\cal A}_2$ be associative superalgebras with parities $\pi_1$ and $\pi_2\,$,
respectively.
Then ${\cal A}={\cal A}_1\otimes {\cal A}_2$
has a natural parity $\pi$ defined by the formula  $\pi(a\otimes b)=\pi_1(a)+\pi_2(b)$.

Let $T_i$ be traces on  ${\cal A}_i$.
Clearly, the function $T$ given by the formula $T(a\otimes b)=T_1(a)T_2(b)$ is a trace on ${\cal A}$.

Let $S_i$ be {\it even} supertraces on  ${\cal A}_i$.
Then the function $S$ such that $S(a\otimes b)=S_1(a)S_2(b)$ is an even supertrace on ${\cal A}$.

%{\bf Notation.}
In what follows, we use three types of brackets:
\begin{eqnarray}
[f,g]&=&fg-gf ,   \nonumber \\
\{f,g\}&=&fg+gf ,\nonumber \\
\lbrack f,g\rbrack _{\varkappa} &=& fg -\varkappa^{\pi(f)\pi(g)}gf.\nonumber
\end{eqnarray}

\subsubsection{The superalgebra of observables}
The superalgebra $H_{W(\R)}(\nu)$ of observables
of the rational Calogero model based on the root system $\R$
is a deformation of the skew product%
\footnote{Let ${\cal A}$ and ${\cal B}$ be the superagebras, and ${\cal A}$ is a ${\cal B}$-module.
We say that the superalgebra ${\cal A}\ast {\cal B}$ is a {\it skew product} if
${\cal A}\ast {\cal B}={\cal A} \otimes {\cal B}$ as a superspace and
$(a_1\otimes b_1)\cdot (a_2\otimes b_2)=a_1 b_1(a_2)\otimes b_1 b_2$.}
 of the Weyl algebra and the group algebra
of a finite group
generated by reflections.

We will define it by Definition \ref{defpage}; now let us describe  the necessary ingredients.

%\begin{definition}\label{reflec}
Let $V={\open R}^N$ be endowed with
a non-degenerate symmetric bilinear form $(\cdot,\cdot)$
and the vectors $\vec a_i$ constitute an orthonormal basis in $V$,
i.e.
$$(\vec a_i,\, \vec a_j)=\delta_{ij}.$$
Let
$x^i$ be
the coordinates of $\x\in V$, i.e.  $\x=\vec a_i\,x^i$.
Then $(\x,\,\yy)=\sum_{i=1}^N x^i y^i$ for any $\x,\,\yy \in V$.
The indices $i$ are raised and lowered by means of the forms $\delta_{ij}$
and $\delta^{ij}$.

For any nonzero $\vv \in V={\open R}^N$, define the {\it reflections}
$R_\vv$ as follows:
\be\label{ref}
R_\vv (\x)=\x -2 \frac {(\x,\,\vv)} {(\vv,\,\vv)} \vv \qquad
\mbox{ for any }\x \in V.
\ee

The reflections (\ref{ref}) have the following properties
\be\label{prop}
R_\vv (\vv)=-\vv,\qquad R_\vv^2 =1,\qquad
({R}_\vv (\x),\,\vec u)=(\x,\,{R}_\vv (\vec u))\quad
\mbox{ for any }\vv,\,\x,\,\vec u\in V.
\ee

%\begin{definition}\label{rootsys}

A finite set of vectors $\R\subset V$ is said to be a {\it root system}
if the following conditions hold:

i) $\R$ is ${R}_\vv$-invariant for any $\vv \in \R$,

ii) if $\vv_1,\vv_2\in \R$ are collinear, then either $\vv_1=\vv_2$ or $\vv_1=-\vv_2$.
%\end{definition}

\medskip
The group $W(\R)\subset O(N, {\open R})\subset End(V)$ generated by all reflections ${R}_\vv$ with
$\vv \in \R$ is finite.

%%%%%%%%%%%%%%%%%%%%%%%%%%%
{
\medskip
As it follows from this definition of a root system, we consider
both crystallographic and non-crystallographic root systems.
We consider also the empty root system denoted by $A_0$,
assuming that it generates the trivial group consisting of the unity element only.
}
%%%%%%%%%%%%%%%%%%%%%%%%%%%%

\medskip

Let ${\cal H}^\alpha$, where $\alpha=0,1$\,, be two copies of $V$ with
orthonormal bases $a_{\alpha\,i} \in {\cal H}^\alpha$, where $i=1,\,...\,,\,N$.
For every vector $\vv=\sum_{i=1}^N {\vec a_i} v^i \in V$, let
$v_\alpha\in {\cal H}^\alpha$
be the vectors $v_\alpha=\sum_{i=1}^N  a_{\alpha \, i} v^i$,
so four bilinear forms on ${\cal H}^0\oplus {\cal H}^1$
can be defined by the expression
\be
\label{bill}
(x_\alpha,\,y_\beta)=(\vec x,\, \vec y) \mbox{ for } \alpha,\beta = 0,1\,,
\ee
where $\vec x,\,\vec y \in V$ and $x_\alpha,\,y_\alpha \in {\cal
H}^\alpha$ are their copies.
The reflections $R_\vv$ act on ${\cal H}^\alpha$ as follows:
\be\label{refl}
R_\vv(h_\alpha)=h_\alpha -
2\frac{(h_\alpha,\,v_\alpha)}{(\vv,\,\vv)}v_\alpha
\qquad \mbox{ for any } h_\alpha\in {\cal H}^\alpha.
\ee
So the $W(\R)$-action on the spaces ${\cal H}^\alpha$
is defined.

\medskip
Let ${\open C}[W(\R)]$ be the group algebra of $W(\R)$, i.e., the
set of all linear combinations $\sum_{g\in W(\R)} \alpha_g \bar g$,
where $\alpha_g \in {\open C}$ and we temporarily use the notations $\bar g$
to distinguish $g$ considered as an element of $W(\R)\subset End(V)$
from the same element $\bar g \in {\open C}[W(\R)]$
of the group considered as an element of the group algebra.
The addition in ${\open C}[W(\R)]$ is defined as follows:
$$
\sum_{g\in W(\R)} \alpha_g \bar g + \sum_{g\in W(\R)} \beta_g \bar g
= \sum_{g\in W(\R)} (\alpha_g + \beta_g) \bar g
$$
and the multiplication is defined by setting
$\overline {g_1\!}\,\, \overline {g_2\!} = \overline {g_1 g_2}$.

Note that the additions in ${\open C}[W(\R)]$ and in $End (V)$ differ.
For example, if $I\in W(\R)$ is unity and the matrix $K=-I$ from $End(V)$
belongs to $W(\R)$, then $I+K=0$ in $End(V)$ while
$\overline I + \overline K \ne 0$ in ${\open C}[W(\R)]$.
In what follows, the element $\overline K\in H_{W(\R)}(\nu)$ is
a {\it Klein operator}.
\footnote{Let ${\cal A}$ be an associative superalgebra with parity $\pi$.
Following M.Vasiliev, see, e.g. \cite{V},
we say that an element $K\in {\cal A}$ is a {\it Klein operator} if $\pi(K)=0$,
$Kf=(-1)^{\pi(f)}fK$ for any $f\in {\cal A}$ and $K^2=1$. Every Klein operator
belongs to the {\it anticenter} of the superalgebra ${\cal A}$, see \cite{sosrus}, p.41.

Any Klein operator, if exists, establishes an isomorphism between the space of
even traces and the space of even supertraces on ${\cal A}$.
Namely, if $f\mapsto T(f)$ is an even trace, then $f\mapsto T(fK)$
is a supertrace, and if $f\mapsto S(f)$ is an even supertrace, then $f\mapsto S(fK)$ is a trace.

It is proved in \cite{all} that if $H_{W(\R)}(\nu)$ has isomorphic spaces of the traces and supertraces,
then $H_{W(\R)}(\nu)$ contains a Klein operator.  }

\medskip
Let $\nu$ be a set of constants $\nu_\vv$ with $\vv\in\R$ such that
$\nu_\vv=\nu_{\vec w}$ if $R_\vv$ and $R_{\vec w}$ belong to
one conjugacy class of $W(\R)$.

 \medskip

\begin{definition}%\label{HW}
%{\bf Definition.}
\label{defpage}
The superalgebra $H_{W({\R})}(\nu)$ is an associative superalgebra with unity ${\mathbf 1}$
of polynomials in the $a_{\alpha\,i}$ with coefficients
in the group algebra ${\open C}[W(\R)]$ subject to the relations
\bee
\overline g h_\alpha&=&g(h_\alpha)h_\alpha \overline g  
%&\ &
\mbox{ for any } g\in W(\R)
                   \mbox{ and } h_\alpha \in {\cal H}^\alpha \\
\label{rel}
 \!\!\!\!\!\! [ x_\alpha \overline I  , y_\beta \overline I ] &=& \varepsilon_{\alpha\beta}
       \left((\vec x,\, \vec y) {\bar 1}  \overline I +
       \sum_{\vv\in\R} \nu_\vv
\frac {(\vec x,\,\vv)(\vec y,\,\vv)}{(\vv,\,\vv)}{\bar 1} \overline {R_\vv} \right)
% &\ &
\mbox{ for any  $x_\alpha \in {\cal H}^\alpha$ and $y_\beta\, \in {\cal H}^\beta$}.
\eee
where $\varepsilon_{\alpha\beta}$
is the antisymmetric tensor, $\varepsilon_{01}=1$, and $\bar 1$ is the unity in ${\open C}[a_{\alpha\,i}]$.
The element   ${\bf 1}=\bar 1 \cdot \overline I$ is the unity  of $H_{W(\R)}(\nu)$.% 
\footnote{Clearly, $H_{W(\R)}$ does not contain  either $\bar 1\in {\open C}$ or $\overline I$.}
The action of any operator $g\in End (V)$ is given by a matrix $g_i^j$:
\bee\label{ga}
g(a_{\alpha \,i} h^i)&=&a_{\alpha \,i} g_i^j h^j, \quad g_1(g_2( h_\alpha))=(g_1 g_2)(h_\alpha)
\ \mbox{ for any } h_\alpha= a_{\alpha \,i} h^i \in {\cal H}^\alpha,\\
g({\bar 1})&=&{\bar 1}.                                      \label{g1}
\eee
The commutation relations (\ref{rel}) suggest
to define the {\it parity} $\pi$ by setting:
\be
\pi (a_{\alpha\,i}\overline g)=1
\ \mbox{ for any }\alpha,\ i \mbox{ and }g\in \G
\qquad \pi(\bar 1 \overline g)=0 \ \mbox{ for any } g\in \G.
\ee

We say that $H_{W({R})}(\nu)$ is a {\it the superalgebra of
observables of the Calogero model based on the root system $\R$}.
\end{definition}

\medskip
These algebras (with parity forgotten) are particular cases of {\it Symplectic Reflection Algebras}
\cite{sra} and are also known as {\it rational Cherednik algebras}.

\medskip
Below we will usually designate ${\mathbf 1}$, $\bar 1$, $I$ and
$\overline  I$ by $1$, and $F\overline I=\bar I F$ by $F$ for any $F\in {\open C}[a_{\alpha\,i}]$,
and ${\bar 1} G=G{\bar 1}$ by $G$ for any $G\in \G$.
Besides, we will just write $g$ instead of $\overline g$
because it will always be clear, whether $g\in W(\R)$ or $g\in {\open C}[W(\R)]$.

\medskip
The associative algebra $H_{W({R})}(\nu)$ has a faithful representation
via Dunkl differential-difference operators $D_i$, see \cite{Dunkl}, acting on
the space of smooth functions on $V$. Namely, let $v_i=\delta_{ij}v^j$, $x_i=\delta_{ij}x^j$,
\be\label{Dun}
D_i=
\frac {\partial} {\partial x^i} + \frac{1}{2}  \sum_{\vv\in\R} \nu_\vv\frac {v_i} {(\x,\,\vv)} (1-R_\vv)
\ee
and
\cite{Poly, BHV}
\be\label{aa}
a_{\alpha\, i} =\frac 1 {\sqrt{2}} (x_i + (-1)^\alpha D_i) \quad \mbox{for }
\alpha =0,1.
\ee

The reflections $R_\vv$ transform the deformed creation and annihilation
operators
(\ref{aa})
as vectors:
\be\label{comav}
R_\vv a_{\alpha\,i} = \sum_{j=1}^N \left(\delta_{ij} - 2
\frac {v_i v_j}{(\vv,\,\vv)}\right)a_{\alpha\,j}  R_\vv.
\ee
Since $[D_i,\, D_j]=0$, see \cite{Dunkl}, it follows that
\be\label{comaa}
[a_{\alpha\, i}, a_{\beta\, j}] = \varepsilon_{\alpha\beta}
\left(\delta_{ij}+
\sum_{\vv\in\R} \nu_\vv \frac {v_i v_j}{(\vv,\,\vv)}R_\vv\right),
\ee
which manifestly coincides with (\ref{rel}).

Observe an important property of superalgebra
$H_{W({R})}(\nu)$: the Lie (super)algebra of its inner
derivations
contains  ${\fr sl}_2$ generated by the operators
\be\label{sl2}
T_{\alpha\beta}= \frac{1}{2}  \sum_{i=1}^N
\left\{a_{\alpha\, i},\,a_{\beta\, i}\right\}
\ee
which commute with ${\open C}[W(\R)]$, i.e.,
$[T_{\alpha\beta},\,R_\vv]=0$,
and act on $a_{\alpha\, i}$ as on vectors
of the irreducible 2-dimensional
${\fr sl}_2$-modules:
\be\label{sl2vec}
\left[T_{\alpha\beta},\,a_{\gamma\, i}\right]=
\varepsilon_{\alpha\gamma}a_{\beta \, i} +
\varepsilon_{\beta\gamma}a_{\alpha\, i}, \quad \mbox{ where } i=1,\,\dots\,, N.
\ee

The restriction of the operator
$T_{01}$ in the representation (\ref{aa}) on the subspace
of $W(\R)$-invariant functions on $V$
is a second-degree differential operator which is the well-known
Hamiltonian of the rational Calogero model, see \cite{Cal}, based on the
root system $\R$, see \cite{OP}.
One of the relations (\ref{sl2}), namely,
$[T_{01},\,a_{\alpha\, i}]=
-(-1)^\alpha a_{\alpha\, i}$, allows one to find the solutions
of
the equation $T_{01}\psi =\epsilon \psi$ and eigenvalues $\epsilon$
via usual Fock procedure with
the vacuum $ |0\rangle$ such that
$a_{0\,i}|0\rangle$=0 $\mbox{ for any } i$, see  \cite{BHV}.
After $W(\R)$-symmetrization these eigenfunctions
become the wave functions of the Calogero Hamiltonian.

%%%%%%%%%%%%%%%%%%%%%%%%%%%%%%%%%%%%%%%%%%%%%%%%%%%%%%%%%%%%%%%%%

\section{The $\varkappa$-traces on $H_{W({R})}(\nu)$}\label{trace}

Every $\varkappa$-trace $\str (\cdot)$ on ${\cal A}$ generates
the following bilinear form on ${\cal A}$:
\bee\label{bf}
B_{\str }(f,g)=\str (f\cdot g) \mbox{ for any } f,g\in {\cal A}.
\eee

It is obvious that if such a bilinear form $B_{\str }$ is degenerate,
then the
null-vectors of this form (i.e., $v \in {\cal A}$ such that $B(v, x)=0$ for any
$x\in {\cal A}$)
constitute the two-sided
ideal ${\cal I}\subset {\cal A}$.
If the $\varkappa$-trace generating degenerate bilinear form is homogeneous (even or odd),
then the corresponding ideal is a superalgebra.

If $\varkappa=-1$, the ideals of this sort are present, for example,
in the superalgebras $H_{W(A_1)}(\nu)$
(corresponding to the two-particle Calogero model) at $\nu =k+ \frac{1}{2} $,
see \cite{V}, and in the superalgebras $H_{W(A_2)}(\nu)$
(corresponding to three-particle Calogero model)
at $\nu =k+ \frac{1}{2} $
and $\nu=k\pm\frac1 3$, see \cite{K2},
for every integer $k$. For all
other values of $\nu$ all supertraces on these superalgebras
generate nondegenerate bilinear forms (\ref{bf}).

%%%%%%%%%%%%%%
{
The general case of $H_{W(A_{n-1})}(\nu)$ for arbitrary $n$ is considered
in \cite{IL}. Theorem 5.8.1 of \cite{IL} states that
the associative algebra $H_{W(A_{n-1})}(\nu)$ is not simple if and only if $\nu=\frac{q}{m}$, where $q,m$
are mutually prime integers, and $1<m\leqslant n$, and presents the structure of corresponding ideals.

%%%%%%%%%%%%%%%%%%
Conjecture: Each of the ideals found in \cite{IL} is
the set of null-vectors of the degenerate bilinear form (\ref{bf})
for some $\varkappa$-trace $sp$
on $H_{W(A_{n-1})}(\nu)$.%
}%
\footnote
{
{The dimension of the space of supertraces on $H_{W(A_{n-1})}(\nu)$ is the number
of the partition of $n\geqslant 1$ into the sum of different positive integers, see \cite{KV},
and the space of the traces on $H_{W(A_{n-1})}(\nu)$ is one-dimensional
for $n\geqslant 2$ due to Theorem \ref{main1}, see also \cite{all}.
}}
%%%%%%%%%%%%%%%%%%

\subsection{Main results}

\indent
\begin{theorem}
{\it Each nonzero $\varkappa$-trace on $H_{W(\R)}(\nu)$ is even.}
\end{theorem}

\indent\begin{proof}
The space of superalgebra $H_{W(\R)}(\nu)$ can be decomposed into
the direct sum of irreducible
${\fr sl}_2$-modules (Lie algebra ${\fr sl}_2$ is defined by eq. (\ref{sl2})).
Clearly, each $\varkappa$-trace should vanish on all these irreducible
modules
except singlets, and can take nonzero value only on singlets, i.e.,
on elements $f\in H_{W(\R)}(\nu)$ such that $[T_{\alpha\beta},\, f]=0$ for $\alpha,\beta=0,1$.
So, if $\str(f)\ne 0$, then $[T_{0\,1},\, f]=0$, which implies $\pi(f)=0$.
\end{proof}

\begin{theorem}\label{main0}
{\it
The dimension of the space of $\varkappa$-traces on the superalgebra $H_{W(\R)}(\nu)$
is equal to
the number of conjugacy
classes of elements without eigenvalue $\varkappa$ belonging to the Coxeter
group $W(\R)\subset End({\open R}^N)$ generated by the finite root system $\R\subset{\open R}^N$.
}
\end{theorem}

\indent\begin{proof}
This Theorem follows from Theorem \ref{th6} and Theorem \ref{th5}.
\end{proof}

Clearly, Theorem \ref{main0} implies the following theorem

\begin{theorem}\label{main1}
{\it
Let the Coxeter
group $W(\R)\subset End({\open R}^N)$ generated by the finite root system $\R\subset{\open R}^N$
have $T_\R$ conjugacy classes without eigenvalue $1$ and $ST_\R$
conjugacy classes without eigenvalue $-1$.

Then the superalgebra
 $H_{W(\R)}(\nu)$ possesses
$T_\R$ independent traces and $ST_\R$ independent supertraces.
}
\end{theorem}

\section{Ground Level Conditions}\label{anal1}

Clearly, the $\bar 1 \cdot {\open C}[W(\R)]$
is a subalgebra of $H_{W({R})}(\nu)$ isomorphic to $\G$.

It is easy to describe all $\varkappa$-traces on
${\open C}[W(\R)]$. Every $\varkappa$-trace on ${\open C}[W(\R)]$
is completely determined by its values on
$W(\R)
%\subset {\open C}[W(\R)]
$
and  is a central function on $W(\R)$, i.e.,
the function constant on the conjugacy classes
due to $W(\R)$-invariance.
Thus, the number of the $\varkappa$-traces on $\G$
is equal to the number of conjugacy classes in $W(\R)$.

Since $\G \subset H_{W(\R)}(\nu)$, some
additional restrictions on these functions
follow from the definition (\ref{scom}) of $\varkappa$-trace and the defining relations (\ref{rel}) for
$H_{W(\R)}(\nu )$.  Namely, consider $g\in {W(\R)}$
and elements $c_i^\alpha\in {\cal H}^\alpha$ such that
\be\label{eigss}
\sigmama c_i^\alpha=\varkappa c_i^\alpha \sigmama.
\ee
Then, eqs. (\ref{scom}) and (\ref{eigss}) imply that
$$
\str    \left ( c_i^0 c_j^1 \sigmama
\right )=  \varkappa  \str    \left ( c_j^1 \sigmama c_i^0\right )=  \str    \left ( c_j^1 c_i^0
\sigmama \right ),
$$
and therefore
\be\label{mm}
\str    \left ( [ c_i^0, c_j^1] \sigmama \right )=0.
\ee

Since $[ c_i^0, c_j^1] \sigmama  \in \G $,
the conditions (\ref{mm}) selects the central functions on $\G$ which can
in principle be extended to $\varkappa$-traces on $H_{W(\R)}(\nu)$,
and Theorem \ref{th6} states that each central functions on $\G$,
which satisfy conditions (\ref{mm}), can be
extended to $\varkappa$-trace on $H_{W(\R)}(\nu)$.
In \cite{KV}, the conditions (\ref{mm})  are called {\it Ground Level Conditions}.

Ground Level Conditions (\ref{mm}) is an overdetermined system
of linear equations for the central functions on $\G$.
The dimension of the space of its solution is given in Theorem \ref{th5}.

\subsection{The number of independent solutions of Ground Level Conditions}
\label{anal2}

Let
us introduce the gradation $E$ on the vector space of
${\open C}[W(\R)]$.
For any $ g\in W(\R)$,
consider the subspaces ${\cal E}^\alpha (g) \subset {\cal H}^\alpha$:
\be
\label{eigs}
{\cal E}^\alpha (g) =\{h\in {\cal H}^\alpha\mid g (h)=\varkappa h \}.
\ee
Clearly, $\dim~{\cal E}^0(g)=\dim~{\cal E}^1(g)$.
Set%
\footnote{
It follows from Lemma \ref{l2} formulated below that if $\varkappa=-1$, then $\rho(g)=E(g)|_{mod 2}$
is a parity on the group algebra ${\open C}[W(\R)]$.
It is a well known parity of elements of the Coxeter group $W(\R)$.
Besides $(E(g)|_{\varkappa=+1}-E(g)|_{\varkappa=-1})|_{mod 2}=N|_{mod 2}$\,.
}
\be\label{grad}
E(g)=\dim~\,{\cal E}^\alpha(g).
\ee
For any $g\in W(\R)$, $E(g)$ is an integer such that $0\leqslant E(g) \leqslant N$.

{\bf Notation.} Let $W_l$ denote a subset of all elements of the group $g\in W(\R)$
 such that
$E(g)=l$.

Clearly,
\be\label{wdeco}
W(\R)=\bigcup_{l=0}^N W_l.
\ee

The set $W_l$ is $W(\R)$-invariant,
and we can introduce $W_l^\ast$ -- the space of $W(\R)$-invariant functions
on $W_l$.

\begin{theorem}\label{th5a}
{\it Each function $S\in W_0^\ast$ can be extended uniquely to the
central function on $W(\R)$ satisfying the Ground Level Conditions.}
\end{theorem}

The following theorem follows from Theorem \ref{th5a}:

\begin{theorem}\label{th5}
{\it The dimension of the space of solutions of
Ground Level Conditions (\ref{mm}) is equal to the
number of conjugacy classes in $W(\R)$ with $E(g)=0$.}
\end{theorem}

Theorems \ref{th5a} and \ref{th5} are proved below
simultaneously.

\medskip
The following lemmas are needed to prove these theorems.

\begin{lemma} \label{lemma2}
{\it
Let $g$ be an orthogonal $N\times N$ real matrix
without eigenvalue $\varkappa$,
i.e., the matrix $g-\varkappa$ is invertible.
Then the matrix $R_{\vv}g$ has exactly one
eigenvalue equal to $\varkappa$.
}
\end{lemma}

\indent\begin{proof}
Consider the equation
$R_{\vv}g\x-\varkappa\x=0$ or, equivalently, $g\x-\varkappa R_{\vv}\x=0$
for eigenvector $\x$ corresponding to eigenvalue $\varkappa$.
Using the definition of $R_{\vv}$
this equation can be expressed as
$$
g\x - \varkappa(\x-2\frac {(\vv,\, \x)} {|\vv|^2} v) =0;
$$
hence,
\be\label{x}
\x=-2\varkappa \frac {(\vv,\, \x)} {|\vv|^2}(g-\varkappa)^{-1} \vv.
\ee
It remains to show that this equation has a nonzero solution.
Let $\vv=(g-\varkappa)\vec w$, and it follows from eq. (\ref{x})
that $\x=\mu \vec w$, where $\mu\in {\open R}$.
Then 
\bee
|\vv|^2&=&2(|\vec w|^2-\varkappa(\vec w,\, g\vec w))\,,
\nn
-2\varkappa (\vv, \x) &=&
 2(|\vec w|^2   -\varkappa(\vec w,\, g\vec w))\mu,\nonumber
\eee
and eq. (\ref{x}) becomes an identity $\mu \vec w=\mu \vec w$.
So the vector $\x_1=(g -\varkappa)^{-1} \vv$ is the only solution,
up to a factor.
\end{proof}

\begin{lemma}\label{l2}
{\it
Let $g$ be an orthogonal $N\times N$ real matrix
and $\cc_i$,
where $i=1,\, ... ,\, E(g)$, the complete orthonormal set of its
eigenvectors corresponding to eigenvalue $ \varkappa$.
Then
\\
\indent{\em i)}
$E(R_{\vv} g) = E(g)+1$ if $(\vv,\,\cc_i)=0$ for all $i$;
\\
\indent{\em ii)}
if there exists an $i$ such that $(\vv,\, \cc_i) \neq 0$, then
$E(R_{\vv} g)=E(g)-1$ and the space of the
eigenvectors of $R_{\vv} g$ corresponding to
eigenvalue $\varkappa$ is the subspace of
$span\{\cc_1,\,...,\, \cc_{E(g)}\}$ orthogonal to $\vv$.}
\end{lemma}

\indent\begin{proof}
Let $C  \stackrel{def}{=}  span\{\cc_1,\,...,\, \cc_{E(g)}\}$ and let
  $V=C\oplus B$ be orthogonal direct sum.
Clearly, $gB=B$.

Let us seek null-vector $\vec z$ of the operator $R_{\vv} g -\varkappa$, i.e.,
the solution of the equation
\be\label{rg}
R_{\vv} g \vec z - \varkappa \vec z=0,
\ee
in the form
$\vec z=\cc+\vec b$, where $\cc \in C$ and $\vec b \in B$.
The definition of $R_{\vv}$ and (\ref{rg}) yield
\be\label{rg2}
- \frac 2 {(\vv,\,\vv)}\left(\varkappa (\cc,\,\vv) + (g \vec b,\, \vv) \right)\vv + (g-\varkappa)\vec b=0.
\ee

Represent $\vv$ in the form $\vv = \vv_c+\vv_b$, where
$\vv_c\in C$, $\vv_b \in B$. Let $\vv_b = (g-\varkappa)\vec w$. Then eq.
(\ref{rg}) is equivalent to the system
\bee\label{rg41}
&&- \frac 2 {(\vv,\,\vv)}\left(\varkappa (\cc,\,\vv_c) + (g \vec b,\, (g-\varkappa)\vec w) \right)\vv_c = 0,
\\
\label{rg42}
&&- \frac 2 {(\vv,\,\vv)}\left(\varkappa (\cc,\,\vv_c) + (g \vec b,\, (g-\varkappa)\vec w) \right)\vec w + \vec b=0.
\eee

Consider the two cases:

{\it i)} Let $(\vv, \, \cc_i)=0$ for all $i=1,...,E(g)$. So, $\vv_c=0$, and hence $\vv \in B$.
Then (\ref{rg42}) acquires the form
\be\label{rg3}
- \frac 2 {(\vv,\,\vv)} (g \vec b,\, 
(g-\varkappa)\vec w) \vec w + \vec b=0.
\ee
It is easy to check that $\vec b =\vec w$ is the only nonzero solution of
(\ref{rg3}) orthogonal to $C$.

So, all the solutions of eq. (\ref{rg}) are linear combinations of
the vectors
$\vec z_i=\cc_i$, where $i=1,\,...,\,E(g)$, and $\vec z_{E(g)+1}=\vec w$.

 {\it ii)} Let $\vv_c \ne 0$. Then eq. (\ref{rg41}) gives
\be\label{rg43}
\varkappa (\cc,\,\vv_c) + (g \vec b,\, (g-\varkappa)\vec w)  = 0
\ee
which reduces eq. (\ref{rg42}) to $\vec b=0$ which, in its turn, reduces eq. (\ref{rg43})
to
$
(\cc,\,\vv)=0.
$
\end{proof}

Let ${\cal P}$ be the projection
${\open C}[W(\R)] \rightarrow {\open C}[W(\R)]$
defined as
\be \label{(*)}
{\cal P}(\sum_i \alpha_i
\sigmama _i)=\sum_{i:\, g_i \neq {\bf 1}} \alpha_i \sigmama _i
\quad
\mbox{ for any $g_i\in W(\R)$, $\alpha_i\in {\open C}$}.
\ee

\begin{lemma}\label{lemma4}
{\it
Let $g\in W(\R)$.
Let
$c^\alpha_1, c^\alpha_2\in {\cal E}^\alpha (g)\subset H_{W(\R)}(\nu)$
(i.e., $\sigmama  c^\alpha_1=\varkappa c^\alpha_1 \sigmama $,
$\sigmama c^\alpha_2= \varkappa c^\alpha_2 \sigmama $).
Then
\bee\label{main}
E({\cal P}([c^\alpha_1,\,c^\beta_2]) \sigmama )=E(\sigmama )-1 \ \
\mbox{ for any } g\in W(\R).
\eee
}
\end{lemma}

\indent\begin{proof}
Proof easily follows from the formula
\be
{\cal P}([c^\alpha_1,\,c^\beta_2])=\varepsilon^{\alpha\beta}
\sum_{\vv\in\R} \nu_\vv \frac {(\vec c_1,\,\vv)
(\vec c_2,\,\vv)}{(\vv,\,\vv)}R_\vv\,.
\ee
Indeed, if $(\vec c_1,\,\vv)(\vec c_2,\,\vv)\neq 0$, then
Lemma \ref{l2} implies that $E(R_\vv g)=E(g)-1$.
\end{proof}

\subsection{Proof of Theorems \ref{th5a} and \ref{th5}}
\label{anal3}

Due to Lemma \ref{lemma4} some of the Ground Level Conditions express the
$\varkappa$-trace of elements $g$ with $E(g)=l$
via the $\varkappa$-traces of elements $R_{\vv}g$ with $E(R_{\vv}g)=l-1$:
\bee\label{GLC2}
\str   (g)=-\str (([c^0_i,\, c^1_i]-1)g)
\mbox{ if } (\vec c_i,\,\vec c_i)=1.
\eee

We prove Theorems \ref{th5a} and \ref{th5} using induction on $E(g)$.

The first step is simple: if $E(g)=0$, then $\str (g)$ is an arbitrary
central function.
The next step is also simple: if $E(g)=1$, then there exists
a unique element $c^0_1\in {\cal E}^0(g)$
and a unique element $c^1_1\in {\cal E}^1(g)$ such that
$|c^\alpha_1|=1$ and $gc^\alpha_1=\varkappa c^\alpha_1g$.
Since $(([c^0_1,\, c^1_1]-1)g)\in \G$ and
$E(([c^0_1,\, c^1_1]-1)g)=0$, then
\be\label{sol}
\str (g)=-\str (([c^0_1,\, c^1_1]-1)g)
\ee
is the unique possible value for $\str (g)$ with $E(g)=1$.
In such a way, $W_0^\ast$ is extended to $W_1^\ast$.

A priori these values are not consistent with other Ground Level Conditions.

Suppose that the Ground Level Conditions (\ref{mm})
\be\nonumber
\str  \left ( [ c^0_i, c^1_j] g \right ) =0
\ee
considered for all $g$ with $E(g)\leqslant l$
and for all $c^\alpha_i\in {\cal E}^\alpha (g)$
 such that $(c^\alpha_i,\,c^\beta_j)=\delta_{ij}$,
 where $i=1,,\,...\,,l$,
have $Q_l$ independent solutions.

\begin{statement} \label{state}
{\it The value $Q_l$ does not depend on $l$.}
\end{statement}

\indent\begin{proof}
It was shown above that $Q_1=Q_0$. Let $l\geqslant 1$.
Let us consider $g\in W(\R)$ with $E(g)=l+1$.
Let $c^\alpha_i\in {\cal E}^\alpha(g)$, where $i=1,2$, be such that
$(c^\alpha_i,\,c^\beta_j)=\delta_{ij}$.
These elements $c^\alpha_i$ give the conditions:
\bee
\label{e11} \str (g)=-\str (([c^0_1,\, c^1_1]-1)g), \\
\label{e22} \str (g)=-\str (([c^0_2,\, c^1_2]-1)g), \\
\label{e21} \str ([c^0_1,\, c^1_2]g)=0.
\eee

Below we prove that eqs. (\ref{e11}) and (\ref{e22}) are equivalent
and that eq. (\ref{e21}) follows from them. So, we will prove that eq. (\ref{e11})
considered for all $g\in W_{s}$, where $0 < s \leqslant l+1$
realizes the extension of $W_0^\ast$ to $W_{l+1}^\ast$.

Let us transform (\ref{e11}):
\bee
\str (g)&=&\str (S_1)-\str (S_{12}) \mbox{, where}\\
S_1&=&-\left([c^0_1,\, c^1_1]-1
  -\sum_{\vv\in\R:\,(\vv,\,\vec c_1)(\vv,\,\vec c_2)\neq 0}
\nu_\vv \frac {(\vv,\,\vec c_1)^2}{|\vv|^2}
   R_\vv \right) g =\nn
&{=}&
  -\left(\sum_{\vv\in\R:\,\,(\vv,\,\vec c_1)(\vv,\,\vec c_2)=0}
\nu_\vv \frac {(\vv,\,\vec c_1)^2}{|\vv|^2}
   R_\vv \right) g =\nn
\label{s1}
&{=}&
  -\left(\sum_{\vv\in\R:\,\,(\vv,\,\vec c_2)=0}\nu_\vv
\frac {(\vv,\,\vec c_1)^2}{|\vv|^2}
   R_\vv \right) g \\
S_{12}&=& \left(
   \sum_{\vv\in\R:\,(\vv,\,\vec c_1)(\vv,\,\vec c_2)\neq 0}
\nu_\vv \frac {(\vv,\,\vec c_1)^2}{|\vv|^2}R_\vv
   \right)g.
\eee
It is clear from eq. (\ref{s1}) and Lemma \ref{l2} that $E(S_1)=l$ and
$S_1c^0_2=\varkappa c^0_2S_1$.
Hence, due to eq. (\ref{GLC2}) and inductive hypothesis
\be
\str (S_1)=-\str (([c^0_2,\,c_2^1]-1)S_1)=
\str (([c^0_2,\,c_2^1]-1) (([c^0_1,\,c_1^1]-1)g - S_{12}))
\ee
and as a result
\be
\str (S_1)=
\str (([c^0_2,\,c_2^1]-1)([c^0_1,\,c_1^1]-1)g )
-\str (([c^0_2,\,c_2^1]) S_{12})
+\str (S_{12}).
\ee
Finally, eq. (\ref{e11}) is equivalent under inductive hypothesis to
\be\label{e11a}
\str (g)=
\str (([c^0_2,\,c_2^1]-1)([c^0_1,\,c_1^1]-1)g )
-\str (([c^0_2,\,c_2^1]) S_{12}).
\ee
Analogously, eq. (\ref{e22}) is equivalent under inductive hypothesis to
\be\label{e22a}
\str (g)=
\str (([c^0_1,\,c_1^1]-1)([c^0_2,\,c_2^1]-1)g )
-\str (([c^0_1,\,c_1^1]) S_{21}),
\ee
where
\be
S_{21} = \left(
   \sum_{\vv\in\R:\,(\vv,\,\vec c_1)(\vv,\,\vec c_2)\neq 0}
\nu_\vv \frac {(\vv,\,\vec c_2)^2}{|\vv|^2}R_\vv
   \right)g.
\ee
Now, let us compare the corresponding terms in eqs. (\ref{e11a}) and
(\ref{e22a}).
First, the relation
\be
\str (([c^0_1,\,c_1^1]-1)([c^0_2,\,c_2^1]-1)g )=
\str (([c^0_2,\,c_2^1]-1)([c^0_1,\,c_1^1]-1)g )
\ee
is identically true for every $\varkappa$-trace on ${\open C}[W(\R)]$ since
$[c^0_1,\,c_1^1]$ commutes with $g$.
Second,
\be
\str (([c^0_1,\,c_1^1]) S_{21})=
\str (([c^0_2,\,c_2^1]) S_{12})
\ee
since
\be\label{r12}
\str ( [c^0_1,\,c_1^1] (\vv,\,\vec c_2)^2 R_\vv g)=
\str ( [c^0_2,\,c_2^1] (\vv,\,\vec c_1)^2 R_\vv g)
\ee
for every $\vv\in\R$ such that
$(\vv,\,\vec c_1)(\vv,\,\vec c_2)\neq 0$.
Indeed, the element
\be
\vec c=\alpha \vec c_1 + \beta \vec c_2\,,\mbox{ where }
\alpha=-(\vv,\,\vec c_2)\neq 0 \mbox{ and }
\beta=(\vv,\,\vec c_1)\neq 0\,,
\ee
is orthogonal to $\vv$:
\be\label{ort1}
(\vv,\,\vec c)=0
\ee
and satisfies the relation
\be
R_\vv g c^\alpha = \varkappa  c^\alpha R_\vv g
\ee
due to Lemma \ref{l2}. This fact together with the fact that
\be
E({\cal P}([c^0_i,\,c^1]) R_\vv g) = l-1 \mbox{ for } i=1,2
\ee
(this also follows from Lemma \ref{l2})
and inductive hypothesis imply
\be\label{ort2}
\str ([c^0_i, \,c^1] R_\vv g)=
\str ([c^0, \,c^1_i] R_\vv g)=0 \quad \mbox{ for } i=1,2.
\ee
Substituting
$\vec c_1=\frac 1 \alpha (\vec c-\beta \vec c_2)$
and $\vec c_2=\frac 1 \beta (\vec c-\alpha \vec c_1)$
in the left-hand side of eq. (\ref{r12})
and using eqs. (\ref{ort1}) and (\ref{ort2})
one obtains the right-hand side of eq. (\ref{r12}).
Thus, eq. (\ref{e11}) is equivalent to eq. (\ref{e22});
hence
\be
\str (([c^0_1,\, c^1_1]-1)g)-\str (([c^0_2,\, c^1_2]-1)g)=0
\ee
for every orthonormal pair $c_1,\,c_2\in {\cal E} (g)$.
Consequently,
\be
\str ([c^0_1,\, c^1_2]g)=0
\ee
which finishes the proof of Statement \ref{state} and Theorem \ref{th5}.
\end{proof}

%---------------------------------------------------------

\section{The number of independent $\varkappa$-traces on $H_{W(\R)}(\nu)$}

For proof of the following theorem, see this and subsequent sections.

\begin{theorem} \label{th6}
{\it
Every $\varkappa$-trace on the algebra ${\open C}[W(\R)]$
satisfying the equations
\bee\label{GLC}
\str ([h_0,\,h_1] \sigmama )=0\qquad \mbox{ for any } g\in W(\R) \mbox{
with }E(g)\neq 0 \mbox{ and } h_\alpha \in {\cal E}^\alpha(g),
\eee
can be uniquely extended to a $\varkappa$-trace on $H_{W(\R)}(\nu)$.}
\end{theorem}

\subsection{Notation}

\medskip
For each $g\in W(\R)$, introduce eigenbases $b_{\alpha\,i}$ in ${\open C}\cdot
{\cal H}^\alpha$ ($i=1,...,N$, $\alpha=0,1$)
such that
\bee\label{eig1}
 \sigmama b_{0\, i}=\lambda_i b_{0\, i}  \sigmama  , \\
 \sigmama b_{1\, i}=\frac 1 {\lambda_i} b_{1\, i}  \sigmama  , \label{eig2}\\
(b_{0\, i}, b_{1\, j})=\delta_{ij}\nonumber.
\eee
Let ${\fr B}_g$ be the set of all these $b_{\alpha\, i}$ for a fixed $g$.

In what follows we use the generalized indices $I,J,...$ instead of pairs ($\alpha, i$)
and sometimes write $i(I)$, $\lambda_I$, $\alpha(I)$ meaning that
\be\label{eig3}
b_I = b_{\alpha(I)\, i(I)},\quad  \sigmama b_I = \lambda_I b_I  \sigmama .
\ee

%%%%%%%%%%%%%%%%%%%%%%%%%%%%%%%%%%%%%%%%%%%%%%%%%%%%%%%%%%%%%%%%%%%%%%%%%%%%%%%%%%%%

%%%%%%%%%%%%%%%%%%%%%%%%%%%%%%%%%%%%%%%%%%%%%%%%%%%%%%%%%%%%%%%%%%%%%%%%%%%%%%%%%%%%

Introduce also a symplectic form
\be\label{symp}
{\cal C}_{IJ}=[b_I,\,b_J]|_{\nu=0}
\ee
and let $f_{IJ}$ be the $\nu$-dependent part of the commutator $[b_I,\, b_J]$:
\be\label{fij}
 F_{IJ} \stackrel {def} = [b_I,\,b_J] = {\cal C}_{IJ}+f_{IJ}.
\ee

The indices $I,J$ are raised and lowered with the help of the symplectic
forms $ {\cal C}^{IJ}$ and $ {\cal C}_{IJ}$:
\be\label{rise}
\mu_I=\sum_J{\cal C}_{IJ}\mu^J\,,\qquad \mu^I=\sum_J\mu_J {\cal C}^{JI}\,;
\qquad \sum_M{\cal C}_{IM}{\cal C}^{MJ}=-\delta_I^J\,.
\ee

%%%%%%%%%%%%%%%%%%%%%%%%%%%%%%%%%%%%%%%%%%%%%%%%%%%%%%%%%%%%%%%%%%%%%%%5

Let ${\fr M}(g)$ be the matrix of the map ${\fr B}_{\bf 1}\longrightarrow
{\fr B}_g$
\bee\label{frm}
b_I=\sum_{i,\alpha} {\fr M}^{\alpha\,i}_I(g)\, a_{\alpha\, i}\,.
\eee
Obviously this map is invertible.
Using the matrix notations one can
rewrite (\ref{eig3}) as
\bee\label{eigmat}
\sigmama b_I =\sum_{J=1}^{2N} \Lambda_I^J(g)\,
 b_J \sigmama,
\eee
where the matrix $\Lambda_I^J$ is diagonal, i.e., $\Lambda_I^J=\delta_I^J\lambda_I$.

\medskip
We will say that the monomial $b_{I_1}b_{I_2}\,\dots\,b_{I_k}g$ is {\it regular}
if $b_{I_s}\in {\fr B}_g$ for all $s=1,\dots,k$ and at least one of
$\lambda_{I_s}$ is not equal to $\varkappa$.

\medskip
We will say that the monomial $b_{I_1}b_{I_2}\,\dots\,b_{I_k}g$ is {\it special}
if $b_{I_s}\in {\fr B}_g$ for all $s=1, \dots,k$ and
$\lambda_{I_s}=\varkappa$ for all $s$.
Clearly, that in this case $E(g)>0$.

\medskip
Introduce a lexicographical partial ordering on $\HH$ as follows.
Let $M_1:=P_1(a_{\alpha\, i})g_1$, $M_2:=P_2(a_{\alpha\, i})g_2$,
where $P_1$ and $P_2$ are polynomials, and $g_1,g_2\in W(\R)$.
\be\label{ord}
\text{
We say that $M_1>M_2$ if $\deg P_1>\deg P_2$ or
if $\deg P_1=\deg P_2$ and $E(g_1)>E(g_2)$.
}
\ee

{
\subsection{The $\varkappa$-trace of General Elements}\label{sec5}

To find the $\varkappa$-trace we consider the defining relations (\ref{scom})
as a system of linear equations for the linear function $\str$.

Clearly, this system can be reduced to
\bee\label{ma1}
&& \str \left( [b_I , P (a) \sigmama
]_\varkappa\right)=0\\
\label{ma2}
&&  \str \left(  \tau^{-1} P (a)\sigmama\tau
\right)=
\str \left( P (a) \sigmama
\right)
\eee
for arbitrary polynomials $P$ and arbitrary $g,\tau\in W(\R)$.

Since each $\varkappa$-trace is even, the equation (\ref{ma1}) can be rewritten in the form
\be\label{evenform}
\str \left( b_I  P (a) \sigmama-
\varkappa P (a) \sigmama b_I \right)=0.
\ee

\medskip
{
Clearly, it is possible to express a $\varkappa$-trace of any monomial in $\HH$
in the terms of $\varkappa$-trace on $\G$
using eq. (\ref{evenform}).
Indeed, this can be done in a finite number of the following step operations.

{\bf Regular step operation.}
Let $b_{I_1}b_{I_2}\,\dots\,b_{I_k}g$ be regular monomial, and we may assume without loss
of generality, that $\lambda_{I_1}\ne \varkappa$.

Then
$$
\str(b_{I_1}b_{I_2}\,\dots\,b_{I_k}g)=\varkappa \str(b_{I_2}\,\dots\,b_{I_k}gb_{I_1})=
\varkappa \lambda_{I_1}\str(b_{I_2}\,\dots\,b_{I_k}b_{I_1}g),
$$
which implies
$$
\str(b_{I_1}b_{I_2}\,\dots\,b_{I_k}g)-
\varkappa \lambda_{I_1}\str(b_{I_1}b_{I_2}\,\dots\,b_{I_k}g)
=\varkappa \lambda_{I_1}\str([b_{I_2}\,\dots\,b_{I_k}, \, b_{I_1}]\,g).
$$
Thus,
\be\label{lessdeg}
\str(b_{I_1}b_{I_2}\,\dots\,b_{I_k}g)=
 \frac   {\varkappa \lambda_{I_1}} {1-\varkappa \lambda_{I_1}} \str([b_{I_2}\,\dots\,b_{I_k}, \, b_{I_1}]\,g).
\ee

This step operation expresses the $\varkappa$-trace of any regular
degree $k$ monomial in terms of the $\varkappa$-trace of degree
$k-2$ polynomials.

\medskip
{\bf Special step operation.}
Let $M:=b_{I_1}b_{I_2}\,\dots\,b_{I_k}g$ be special monomial and $E(g)=l>0$.

We can choose a basis $b_I$ in ${\cal E}^0\oplus {\cal E}^1$ such that
${\cal C}_{IJ}|_{{\cal E}^0\oplus {\cal E}^1}$ has the canonical form:
$$
{\cal C}_{IJ}|_{{\cal E}^0\oplus {\cal E}^1}
=\left(
\begin{array}{cc}
0 & I_{E(g)} \\
-I_{E(g)} & 0%
\end{array}%
\right)
$$

Up to a polynomial of lesser degree, the monomial $M$ can be expressed in the form
$$
M=b_{I}^p b_{J}^q \,b_{L_1}\,\dots\,b_{L_{k-p-q}}g + \mbox{ lesser\_degree\_polynomial},
$$
where
\bee
&& 0\leqslant p,q \leqslant k,\quad p+q\leqslant k,\nn
&& \lambda_I=\lambda_J=\lambda_{L_s}=\varkappa \quad \mbox{for any } s, \\
&& {\cal C}_{IJ}=1,\quad {\cal C}_{IL_s}=0,\quad {\cal C}_{JL_s}=0\quad   \mbox{for any } s\,\,. \nonumber
\eee
Let $M':=b_{I}^p b_{J}^q \,b_{L_1}\,\dots\,b_{L_{k-p-q}}\,$.

Now we can derive the equation for $\str (M'g)$. Consider
$$
\str (b_J b_I M'g)=\varkappa \str(b_I M'gb_J)=
\str(b_I M'b_Jg),
$$
which implies
\be\label{so}
\str([b_I M',\,b_J]g)=0.
\ee

Transform the expression $[b_I M',\,b_J]$:
\bee
\!\! \!\! &&\!\! \!\! [b_{I}^{p+1} b_{J}^q \,b_{L_1}\,\dots\,b_{L_{k-p-q}},\, b_J]
=
\nn
\!\! \!\! &&\!\! \!\!
\sum_{t=0}^p b_I^t(1+f_{IJ})b_I^{p-t}b_{J}^q \,b_{L_1}\,\dots\,b_{L_{k-p-q}}
\!\!
+
\!\! \!\!
\sum_{t=1}^{k-p-q}b_{I}^{p+1} b_{J}^q \,b_{L_1}\,\dots\, b_{L_{t-1}}\, f_{L_t\, J}
\, b_{L_{t+1}}\dots\,  b_{L_{k-p-q}}  \,.
\eee
So, eq. (\ref{so}) can be rewritten in the form
\bee\label{000}
\str(M'g)= &-& \str(
\sum_{t=0}^p b_I^tf_{IJ}b_I^{p-t}b_{J}^q \,b_{L_1}\,\dots\,b_{L_{k-p-q}}g   \nn
&+&
\sum_{t=1}^{k-p-q}b_{I}^{p+1} b_{J}^q \,b_{L_1}\,\dots\, b_{L_{t-1}}\, f_{L_t\, J}
\, b_{L_{t+1}}\dots\,  b_{L_{k-p-q}}g),
\eee
which is the desired equation for $ \str(M'g) $.

Due to Lemma \ref{l2} it is easy to see that eq. (\ref{000}) can be rewritten in the form
\be\label{step2}
\str(M'g)=\sum_{\tilde g\in W(\R):\, E(\tilde g)=E(g)-1} \str(P_{\tilde g} (a_{\alpha\,i})\tilde g),
\ee
where the $P_{\tilde g}$ are some polynomials such that $\deg P_{\tilde g} \leqslant \deg M'$.

So, the special step operation expresses the $\varkappa$-trace of a special
polynomial in terms
of the $\varkappa$-trace of polynomials  lesser in the sense of the ordering (\ref{ord}).

Thus, we showed that it is possible to express a  $\varkappa$-trace of any polynomial
in the terms of $\varkappa$-trace on $\G$ using a finite number of regular and special step operations.
Since each step operation is manifestly $W(\R)$-invariant, and the $\varkappa$-trace on $\G$ is
$W(\R)$-invariant also, the resulting $\varkappa$-trace is $W(\R)$-invariant.

This is not a proof of Theorem \ref{th6} yet because the resulting values
of $\varkappa$-traces may
a priori depend on the sequence of step operations used and impose an additional constraints on
the values of $\varkappa$-trace on $\G$.

Below we prove that the value of $\varkappa$-trace does not indeed depend
on the sequence of step operations used. We use the following inductive procedure:

{\bf $(\star)$} Let $F:=P(a_{\alpha\,i})g\in \HH$, where $P$ is a polynomial such that $\deg P=2k$
and $g\in W(\R)$.
Assuming that $\varkappa$-trace is correctly defined for all elements of $\HH$ lesser
than  $F$ relative to the ordering (\ref{ord}), we prove that $\str (F)$ is defined also
without imposing an additional constraints on the solution of the Ground Level Conditions.

The central point of the proof is consistency conditions
(\ref{c1}), (\ref{c2}) and (\ref{c3}) proved in Appendices \ref{appb} and \ref{appc}.
}
}

%%%%%%%%%%%%%%%%%%%%%%%%%%%%%%%%%%%%%%%%%%%%%%%%%%%%%%%%%%%%%%%%%%%                                              %%%%
%%%%%%%%%%%%%%%%%%%%%%%%%%%%%%%%%%%%%%%%%%%%%%%%%%%%%%%%%%%%%%%%%%%                                              %%%%
                                                                                                                 %%%%
   \medskip                                                                                                      %%%%
   Assume that the { Ground Level Conditions} hold.                                                              %%%%
   The proof of Theorem \ref{th6} will be given in a constructive way by                                         %%%%
   the following double induction procedure, equivalent to ($\star$):                                                                     %%%%
                                                                                                                 %%%%
   %\noindent                                                                                                    %%%%

   {\bf (i)} Assume that                                                                                         %%%%
   $$                                                                                                            %%%%
   \str \left([b_I , P_p (a) \sigmama                                                                            %%%%
   ]_\varkappa \right) =0 \ \ \mbox{ for any $P_p (a),$ $\sigmama $ and $I$}                                     %%%%
   $$                                                                                                            %%%%
    provided that $b_I \in {\fr                                                                                  %%%%
   B}_\sigmama$                                                                                                  %%%%
   and                                                                                                           %%%%
   $$                                                                                                            %%%%
   \begin{array}{l}                                                                                              %%%%
   \mbox{$\lambda (I) \neq \varkappa$; $p\leqslant k\,$  or}\\                                                   %%%%
   \mbox{$\lambda (I) =\varkappa$, $E(\sigmama  )\leqslant l$, $p\leqslant k\,$ or}\\                            %%%%
   \mbox{$\lambda (I) =\varkappa$; $p\leqslant k-2$}\,,                                                          %%%%
   \end{array}                                                                                                   %%%%
   $$                                                                                                            %%%%
   where $P_p (a)$ is an arbitrary degree $p$ polynomial in                                                      %%%%
   $a_{\alpha\,i} $ and $p$ is odd. This implies that there exists                                               %%%%
   a unique extension of the $\varkappa$-trace such that the same is true for                                    %%%%
   $l$ replaced with $l+1$.                                                                                      %%%%

                                                                                                                 %%%%
   %\noindent                                                                                                    %%%%
   {\bf (ii)} Assuming that                                                                                      %%%%
   $\str \left( b_I  P_p (a) \sigmama-                                                                           %%%%
   \varkappa P_p (a) \sigmama b_I \right)=0$ for any $P_p (a)$, $\sigmama $ and                                  %%%%
   $b_I \in {\fr                                                                                                 %%%%
   B}_\sigmama$,  where $p\leqslant k$,                                                                          %%%%
   one proves that there exists a unique extension of the $\varkappa$-trace  such                                %%%%
   that                                                                                                          %%%%
   the assumption {\bf (i)} is true for $k$ replaced with $k+2$ and $l=0$.                                       %%%%

                                                                                                                 %%%%
   As a result, this inductive procedure uniquely extends any solution of the                                    %%%%
   { Ground Level Conditions}                                                                                    %%%%
   to a $\varkappa$-trace on the whole $H_{W(\R)} (\nu )$. (Recall                                               %%%%
   that the $\varkappa$-trace of any odd element of $H_{W(\R)} (\nu )$ vanishes                                  %%%%
   because $\varkappa$-trace is even.)                                                                           %%%%

%%%%%%%%%%%%%%%%%%%%%%%%%%%%%%%%%%%%%%%%%%%%%%%%%%%%%%%%%%%%%%%%%%%%%
%%%%%%%%%%%%%%%%%%%%%%%%%%%%%%%%%%%%%%%%%%%%%%%%%%%%%%%%%%%%%%%%%%%%%

%--------------------------------------
\medskip

It is convenient to work with the exponential generating
functions
\be
\label{gf}
\Psi_g   (\mu )=
\str \left ( e^S \sigmama  \right )\,,\mbox{ where }
S= \sum_{L=1}^{2N} (\mu^{L } b_{L} )\,,
\n
\ee
where $g $ is a fixed element of $W(\R)$, $b_L \in {\fr B}_g $ and
$\mu^{L } \in {\open C} $ are independent parameters.  By differentiating eq. (\ref{gf})
with respect to
$\mu^{L }$ one can obtain an arbitrary polynomial in $b_L$ as a coefficient of
$\sigmama $.  The exponential form of the generating functions implies that these
 polynomials are symmetrized.  In these terms, the induction on the degree of
polynomials is equivalent to the induction on the homogeneity degree in
$\mu$ of the power series expansions of $\Psi_g  (\mu )$.

As a consequence of the general properties of the $\varkappa$-trace,
the generating function $\Psi_g  (\mu )$ must be $W(\R)$-invariant:
\bee\label{S_N}
\Psi_{\tau g \tau^{-1}}(\mu)=\Psi_g
(\tilde{\mu})\,,
\eee
where the $W(\R)$-transformed parameters are of the form
\bee
\label{base}
\tilde{\mu}^I=\sum_J  \left({\fr M}(\tau g \tau^{-1})
{\fr M}^{-1}(\tau)\Lambda^{-1}(\tau){\fr M}(\tau)
{\fr M}^{-1}( g )\right)^I_J {\mu}^J
\eee
and matrices ${\fr M}(g )$ and $\Lambda(g)$ are defined in eqs.
(\ref{frm}) and (\ref{eigmat}).

The necessary
and sufficient conditions for the existence of an even $\varkappa$-trace are the
$W(\R)$-covariance conditions (\ref{S_N}) and the condition that
\be\label{start}
\str \left(
\lbrack
b_L  , e^S
\sigmama  \rbrack _\varkappa\right)=0\qquad
\mbox{for any $ g $ and $L$} \,,
\ee
or, equivalently,
\be\label{start1}
\str \left(
b_L  e^ S
\sigmama  - \varkappa e^S
\sigmama b_L \right)=0\qquad
\mbox{for any $ g $ and $L$} \,.
\ee
%%%%%%%%%%%%%%%%%%%%%%%%%%%%%%%%%%%%%%%%%%%%

%%%%%%%%%%%%%%%%%%%%%%%%%%%%%%%%%%%%%%%%%%%%
\subsection{General relations}\label{genrel}

{
{To transform eq. (\ref{start1}) to a form convenient for the proof,}
we use the following
two general relations true for arbitrary operators $X$ and $Y$ and
 parameter $\mu \in {\open C}$:
\be\label{r1}
X \exp(Y+\mu X)=\frac {\partial }{\partial \mu} \exp  (Y+\mu X )+
\int \,t_2 \, \exp (t_1 (Y+\mu X))[X,Y]  \exp (t_2 (Y+\mu X))D^1t ,
\ee
\be\label{r2}
 \exp (Y+\mu X)X=\frac {\partial }{\partial \mu} \exp  (Y+\mu X )-
\int \,t_1 \, \exp (t_1 (Y+\mu X))[X,Y]  \exp (t_2 (Y+\mu X))D^1t
\ee
with the convention that
\be\label{t}
 D^{n-1}t=\delta (t_1 +\ldots +t_n -1)\theta (t_1 )\ldots \theta (t_n )
dt_1 \ldots dt_n \,.
\ee

The relations (\ref{r1}) and (\ref{r2}) can be derived with the help of
partial integration (e.g.,  over $t_1$) and the following formula
\be\label{d}
\frac {\partial }{\partial \mu} \exp  (Y+\mu X ) =
\int \,  \exp (t_1 (Y+\mu X)) X  \exp (t_2 (Y+\mu X))D^1 t\,
\ee
which can be proven by expanding in power series.  The well-known formula
\be
\label{r3}
[X, \exp (Y)]=
\int \,  \exp (t_1 Y)[X,Y]  \exp (t_2 Y)D^1 t
\ee
is a consequence of eqs. (\ref{r1}) and (\ref{r2}).%
\footnote{The independent proof of eq. (\ref{r3}) follows
from the equalities:
$$[X, \exp (Y)]=\lim_{n\rightarrow \infty}[X, (\exp (Y/n))^n]=%
\lim_{n\rightarrow \infty}\sum_{k=0}^{n-1}
(\exp (Y/n))^k [X,(1+\frac 1 n Y)] (\exp (Y/n))^{n-k-1}.
$$
The same trick can be used for the proof of eq. (\ref{d}).
}

With the help of eqs. (\ref{r1}),  (\ref{r2}) and (\ref{eig3}) one rewrites
eq. (\ref{start1}) as
\be\label{nm1}
(1-\varkappa \lambda_L )\frac {\partial }{\partial \mu^L }\Psi_g   (\mu )=
\int
\,(-\varkappa \lambda_L t_1 -t_2 ) \str \ig (
\exp (t_1 S)[b_L ,S]\,\exp (t_2 S)\sigmama \ig )\, D^1 t\,.
\ee
This condition should be true for any $\sigmama $ and $L$ and plays the central
role in the analysis in this section.
Eq. (\ref{nm1}) is an overdetermined system of linear equations
for $\str$; we show below that it has the only solution extending any fixed solution
of the Ground Level Conditions.

There are two essentially distinct cases,  $\lambda_L \neq \varkappa $ and
$\lambda_L =\varkappa $. In the latter case, the eq. (\ref{nm1}) takes the form
\be\label{m1}
0=\int \, \str \ig (
\exp (t_1 S)[b_L ,S]\,exp (t_2 S) \sigmama  \ig )D^1 t\,,\qquad \lambda_L  =\varkappa \,.
\ee

In Appendix \ref{appb} we prove by induction that eqs. (\ref{nm1})
and (\ref{m1}) are consistent in the following sense
\bee\label{c1}
(1 -\varkappa \lambda_K )\frac {\partial }{\partial \mu^K }\int\,(-\varkappa \lambda_L t_1 -t_2 )
\str \ig (
\exp (t_1 S)[b_L ,S]\,\exp (t_2 S) \sigmama
\ig )D^1 t -(L \leftrightarrow K )=0 & & \\
\mbox{for }\lambda_L \neq \varkappa , \ \lambda_K  \neq  \varkappa  & & \nonumber
\eee
and
\be
\label{c2}
(1 -\varkappa \lambda_K )\frac {\partial }{\partial \mu^K }\int \, \str \ig (
\exp (t_1 S)[b_L ,S]\,\exp (t_2 S)\sigmama \ig )D^1 t\,=0
\mbox{ for }   \lambda_L=\varkappa .
\ee
Note that this part of the proof is quite general and does not depend on a
concrete form of the commutation relations between $a_{\alpha\,i}$ in eq.
(\ref{rel}).

By expanding the exponential $e^S$ in eq. (\ref{gf}) into power series in $\mu
^K$ (equivalently $b_K$) we conclude that eq. (\ref{nm1}) uniquely
reconstructs the $\varkappa$-trace of monomials containing $b_K$ with
$\lambda_K\neq \varkappa $
(i.e., {\it regular monomials}) in terms of $\varkappa$-traces of some
lower degree polynomials. Then the consistency conditions (\ref{c1}) and (\ref{c2})
guarantee that eq. (\ref{nm1}) does not impose any additional conditions on
the $\varkappa$-traces of lower degree polynomials and allow one to represent the
generating function in the form
\bee\label{ex1}
\Psi_g  &=& \Phi_g (\mu)\\
 &+&
\sum_{L:\,\lambda_L \neq \varkappa }
\int_0^1
\frac {\mu_L d\tau} {1-\varkappa \lambda_L}\int D^1 t\,(-\varkappa \lambda_L t_1 -t_2 ) \str \ig (
e^{t_1 (\tau S^{\prime\prime}+S^\prime)}[b_L ,(\tau S^{\prime\prime}+S^\prime)]
\,e^{ t_2 (\tau S^{\prime\prime}+S^\prime)}\sigmama \ig )\, ,
\nonumber
\eee
where we
introduced the generating functions $\Phi_g $ for the $\varkappa$-trace of
{\it special polynomials}, i.e.,
the polynomials depending only on $b_L$ with
$\lambda_L=\varkappa $,
\be\label{gff}
\Phi_g   (\mu )\stackrel {def}{=}
\str \left ( e^{ S^\prime} \sigmama  \right ) =
 \Psi_g  (\mu)\ig |_
{(\mu^I=0\ \forall I:\ \lambda_I \neq \varkappa )}
\n
\ee
and
\be
\label{spr}
S^\prime = \sum_{L:\,b_L \in {\fr B}_\sigmama ,\,\lambda_L=\varkappa}
 (\mu^{L } b_L); \qquad S^{\prime\prime}=S-S^\prime\,.
\ee
The relation (\ref{ex1}) successively expresses the $\varkappa$-trace of
higher degree regular
polynomials via the $\varkappa$-traces of lower degree polynomials.

One can see that the arguments above prove the inductive
hypotheses {\bf (i)} and {\bf (ii)} for the particular case where either the
polynomials $P_p (a)$ are regular and/or $\lambda_I \neq \varkappa $.  Note that for
this case the induction {\bf (i)} on the gradation $E$
 is trivial: one simply proves that the degree of the
polynomial can be increased by two.

Let us now turn to a less trivial case of the special polynomials:
\be\label{startprime}
\str \left (b_I   e^{S^\prime} \sigmama - \varkappa e^{S^\prime}\sigmama b_I  \right )=0\,,
\mbox{ where } \lambda_I =\varkappa \,.
\ee
This equation implies
\be\label{startprime1}
\str \left ([b_I,\,   e^{S^\prime}]\sigmama  \right )=0\,,
\mbox{ where } \lambda_I =\varkappa \,.
\ee

Consider the part of $\str \left ([b_I  , \exp S^\prime ] \sigmama  \right )$
which is of degree $k$ in $\mu$ and let $E(g)=l+1$.
By eq. (\ref{m1}) the conditions (\ref{startprime1}) give
\be
\label{m1prime} 0=
\int \,
\str \left( \exp (t_1 S^\prime)[b_I ,S^\prime]\, \exp (t_2 S^\prime) \sigmama
\right) D^1 t\,.  \n
\ee

Substituting $[b_I ,S^\prime]=\mu_I +  \sum_M f_{IM}\mu^M$,
where the quantities $f_{IJ}$ and $\mu_I$ are defined
in eqs. (\ref{fij})-(\ref{rise}), one can rewrite eq. (\ref{m1prime})
in the form
\bee\label{formprime}
\mu_I \Phi_g (\mu ) &=&
- \int
\str \bigg ( \exp (t_1 S^\prime)\sum_M f_{IM}\mu^M\, \exp (t_2 S^\prime) \sigmama
\bigg) D^1 t\,.
\eee

Now we use the inductive hypothesis {\bf (i)}.  The right hand
 side of eq. (\ref{formprime}) is a $\varkappa$-trace of
a  polynomial of degree $\leqslant (k-1)$  in $a_{\alpha\,i}$ in the sector of  degree $k$
polynomials in $\mu$, and $E(f_{IM}g)=l$.
Therefore one can use the inductive hypothesis {\bf
(i)} to obtain the equality
$$
\int \str \ig (  \exp (t_1 S^\prime)\sum_M f_{IM}\mu^M\,  \exp (t_2 S^\prime)
\sigmama  \ig )D^1t =
\int \, \str \ig (   \exp (t_2 S^\prime) \exp (t_1 S^\prime)
\sum_M f_{IM}\mu^M \sigmama  \ig )D^1 t,
$$
where we used
that
$\str (S^\prime F \sigmama )$ $=$ $ \varkappa \str (F
\sigmama  S^\prime) $= $\str (F  S^\prime \sigmama )$ by definition of
$S^\prime$.

As a result, the inductive hypothesis allows one to transform eq. (\ref{startprime})
to the following form:
\be
\label{p9}
X_I  \stackrel{def}{=}  \mu_I \Phi_g (\mu )
+ \str \bigg( \exp (S^\prime ) \sum_Mf_{IM}\mu^M\sigmama  \bigg)=0 \,.
\ee

By differentiating this equation with respect to $\mu^J$ one obtains after
symmetrization
\be
\label{p10}
\frac {\partial }{\partial \mu^J}   \left(
\mu_I
\Phi_g  (\mu )\right)
+(I\leftrightarrow J )=-
\int \str \ig (e^{t_1 S^\prime } b_Je^{t_2 S^\prime }
\sum_M f_{IM}\mu^M \sigmama  \ig )D^1 t +(I\leftrightarrow J ).
\ee

An important point is that the system of equations (\ref{p10}) is equivalent
to the original equations (\ref{p9}) except for the ground level part
$\Phi_g  (0)$.
This can be easily seen from the simple fact that the general solution
of the system of equations
for
entire functions
$X_I(\mu)$
$$
\frac {\partial }{\partial \mu^J} X_I(\mu) + \frac {\partial }{\partial \mu^I} X_{J}(\mu) =0
$$
is of the form
$$X_I(\mu)=X_I(0)+\sum_{J}c_{IJ}\mu^J$$
where
$X_I(0)$ and
$c_{JI}$=$-c_{IJ}$
are some constants.

The part of eq. (\ref{p9}) linear in $\mu$ is however equivalent to the
Ground Level Conditions
analyzed in Section \ref{anal1}.
Thus, eq. (\ref{p10}) contains all information of eq. (\ref{mm})
additional to the Ground Level Conditions.
For this reason we will from now on analyze the
equation (\ref{p10}).

Using again the inductive hypothesis
we move $b_I$ to the left and to
the right of the right hand side of eq. (\ref{p10})
with equal weights equal to $\frac 1 2$ to get
\bee\label{p11}
&{}&
\frac {\partial }{\partial \mu^J}\mu_I \Phi_g  (\mu )+(I\leftrightarrow J )=
-\frac 1 {2} \sum_{M}\str \ig (  \exp (S^\prime )\{b_J ,f_{IM}\}\mu^M \sigmama \ig )\nn
&{}&
-\frac {1}{2}\int \,\sum_{L,M}(t_1 -t_2 ) \str \ig (\exp (t_1 S^\prime )
F_{JL}\mu^L \exp (t_2 S^\prime ) f_{IM}\mu^M \sigmama  \ig )D^1 t
+ (I\leftrightarrow J )   \,.
\eee
The last term in the right hand side of this expression can be shown to
vanish under the $\varkappa$-trace
due to the factor  $t_1 -t_2 $, so
that one is left with the equation
\be\label{dm1}
L_{IJ}\Phi_g  (\mu )=  -\frac {1}{2}
R_{IJ} (\mu )\,,
\ee
where
\bee\label{RIJ}
R_{IJ} (\mu )=\sum_{M}
\str \ig (  \exp (S^\prime )\{b_J ,f_{IM}\}\mu^M \sigmama  \ig )
+(I\leftrightarrow J )
\eee
and
\bee\label{LIJ}
L_{IJ}=\frac {\partial }{\partial \mu^J}\mu_I + \frac {\partial }{\partial \mu^I}\mu_J\,.
\eee

The differential operators $L_{IJ}$ satisfy the standard
$sp(2E(g))$
commutation
relations
\be\label{lcom}
[L_{IJ},L_{KL}]= - \left(
{\cal C}_{IK}L_{JL}+
{\cal C}_{IL}L_{JK}+
{\cal C}_{JK}L_{IL}+
{\cal C}_{JL}L_{IK} \right)
\,.
\ee
In Appendix \ref{appc} we show by induction  that this $sp(2E(g))$ Lie algebra
of differential operators is consistent
with the right-hand side of the basic relation (\ref{dm1}), i.e., that
\be\label{c3}
[L_{IJ},\,R_{KL}]-
[L_{KL},\,R_{IJ}]=    -\left(
{\cal C}_{IK}R_{JL}+
{\cal C}_{JL}R_{IK}+
{\cal C}_{JK}R_{IL}+
{\cal C}_{IL}R_{JK}     \right)
\,.
\ee
}

%%%%%%%%%%%%%%%%%%%%%%%%%%%%%%%%%%%%%%%%%%%%

{
Generally, these consistency conditions guarantee that
eqs. (\ref{dm1})
express $\Phi_g  (\mu ) $ in terms of $R^{IJ}$ in the following
way
\bee
\label{ex2}
\Phi_g (\mu)&=&
\Phi_g (0)+\frac  {1}{8E(g)}\sum_{I,J=1}^{2E(g)} \int_0^1
\frac{dt}{t} (1-t^{2E(\sigmama  )})
(L_{IJ} R^{IJ})(t\mu ) \,,
\eee
provided
\be
\label{0}
R^{IJ}(0)=0\,.
\ee
The latter condition must hold for the consistency of eqs. (\ref{dm1})
since its left hand side vanishes at $\mu^I =0$. In the expression (\ref{ex2})
it guarantees that the integral over $t$  converges.
In the case under consideration the property eq. (\ref{0})
is indeed true as a consequence of the definition (\ref{RIJ}).

Taking into account  Lemma \ref{l2} and the explicit form (\ref{RIJ}) of $R_{IJ}$
 one concludes that eq.
(\ref{ex2}) uniquely expresses the $\varkappa$-trace
of special polynomials in terms of the
$\varkappa$-traces of polynomials of lower degrees or in terms of the $\varkappa$-traces of special
polynomials of the same degree multiplied by elements of $W(\R)$ with a lower value of $E$
provided that the $\mu$-independent term $\Phi_g (0)$ is an arbitrary solution
of {the Ground Level Conditions}.  This completes the proof of Theorem \ref{th6}. $\blacksquare$
}

\section{Undeformed skew product of Weyl superalgebra and finite group generated by
reflections ($H_{W(\R)}(0)$)}

Consider $H_{W(\R)}(0)$. This algebra is the skew product of the Weyl superalgebra
and the group algebra of the finite group $W(\R)$ generated by a root system $\R\subset V={\open R}^N$.
Algebras of this type, and their generalization, were considered in \cite{pass}.

The superalgebra $H_{W({\R})}(0)$ is an associative superalgebra
of polynomials in $a_{\alpha\, i}$\,, where $\alpha=0,1$ and $i=1,...,N$, with coefficients
in the group algebra ${\open C}[W(\R)]$ subject to the relations
\bee\label{rel01}
g a_{\alpha\, i} &=& \sum_{k=1}^N g^k_i a_{\alpha\, k}g  \mbox{ for any } g\in W(\R)
                   \mbox{ and } a_{\alpha\, i} \,,\\
\label{rel02}
  [ a_{\alpha\,i}, a_{\beta\, j}] &=& \varepsilon_{\alpha\beta}
       \delta_{ij} \,,
\eee
where $\varepsilon_{\alpha\beta}$
is the antisymmetric tensor, $\varepsilon_{01}=1$, and $g_i^k$ is a matrix
realizing representation of $g\in W(\R)$ in $End(V)$.
The commutation relations (\ref{rel01})-(\ref{rel02}) suggest
to define the {\it parity} $\pi$ by setting:
\be
\pi (a_{\alpha\, i})=1
\mbox{ and } \pi(g)=0 \ \mbox{ for any } g\in W(\R).
\ee

Unifying indices $i$ and $\alpha$ in one index $I$ one can rewrite eq. (\ref{rel02})
as
\be\label{rel03}
[a_I,\,a_J]=\omega_{IJ}\,,
\ee
where $\omega_{IJ}$ is a symplectic form.

\medskip

It is 
easy to find
the general solution of eqs. (\ref{nm1}) and (\ref{m1}) for generating function
of $\varkappa$-traces:

\begin{enumerate}
\item
If $g\in W(\R)$ and $E(g)\ne 0$, then $sp(P(a_I)g)=0$ for any polynomial $P$.

\item
If $g\in W(\R)$ and $E(g)= 0$, then $sp(g)$ is
an arbitrary central function on $W(\R)$.

\item\label{lamb}
Let $E(g)= 0$. There exists a complete set $b_{\alpha, k}$  of eigenvectors of $g$ for each $\alpha$,
such that $g b_K=\Lambda_ K b_K g$
and ${\cal C}_{KL}=[b_K,\,b_L]$ is nondegenerate
skewsymmetric form such that ${\cal C}_{KL} \ne 0$ only if $\lambda_K \lambda_L=1$.
In this notation, let
\begin{eqnarray*}
S(\mu,b)&=& \sum_K \mu^K b_K, \\
Q(\mu)&=& \frac 1 4 \sum_{KL} \mu^K \mu^L \tilde{\cal C}_{KL}\,,
\end{eqnarray*}
where
\begin{align*}
\tilde{{\cal C}}_{KL}=-\frac{1+\varkappa\lambda_{K}}{1-\varkappa\lambda_{K}}{\cal C}_{KL}
=\tilde{{\cal C}}_{LK}\,.
\end{align*}%

Then
\be\label{ans}
sp\left(e^{S(\mu,b)} g \right) = e^{Q(\mu)} sp(g)\,.
\ee

\end{enumerate}

\medskip

%%%%%%%%%%%%%%%%%%%%%%%%%%%%%%%%%%%%%%%%%%%%%%%%%%%%%%%%%%%%%%%%%%%%%%%%%

The solution eq. (\ref{ans}) can be obtain in initial basis also.

Let $S=\sum_{\alpha\,i} \mu^{\alpha\,i} a_{\alpha\,i}$\,,
$
\Psi(g,\mu,t)=\str( e^{tS}g),\;\Psi(g,\mu)=\str ( e^{S}g )
=\Psi(g,\mu,1).
$
Then
\begin{align}\label{nu0}
\str\left(  [a_{\alpha i},e^{tS}g]_{\varkappa}\right)   &  =\str \left(
t\varepsilon_{\alpha\beta}\delta_{ij}\mu^{\beta j}e^{tS}g+e^{tS}a_{\alpha
j}gp_{i}^{j}\right)  ,\;\mbox{ where } p_{i}^{j}=(1-\varkappa g)_{i}^{j},
\end{align}%

Since $E(g)=0$, the matrix $p_{i}^{j}$ is invertible, so eq.
(\ref{nu0}) gives

\begin{align*}
&  \frac d {dt}{\Psi}(g,\mu,t)=-\mu^{\alpha i}\varepsilon_{\alpha\beta}q_{i}%
^{k}\delta_{kj}\mu^{\beta j}\Psi(g,\mu,t),
 \mbox{ where }
  q_{i}^{k}=\left(  \frac{1}{I-\varkappa g}\right)  _{i}^{k}=\frac{1}%
{2}\left(  \frac{I+\varkappa g}{I-\varkappa g}\right)  _{i}^{k}+\frac{1}%
{2}\delta_{i}^{k}
\end{align*}%

So

\begin{align*}
\frac d {dt}{\Psi}(g,\mu,t)  &  =-\mu^{\alpha i}\varepsilon_{\alpha\beta}%
\tilde{\omega}_{ij}\mu^{\beta j} {\Psi}(g,\mu,t),
\mbox{ where } \tilde{\omega}_{ij}=\frac{1}{2}\left(
\frac{1+\varkappa g}{1-\varkappa g}\right)  _{i}^{k}\delta_{kj}=-\tilde{\omega}_{ji}%
\end{align*}
and finally
\be\nonumber
\Psi(g,\mu)    = \exp
\Big (
-\frac{1}{2}\mu^{\alpha i}\varepsilon_{\alpha\beta}\tilde{\omega}_{ij}\mu^{\beta j}
\Big )
\str (\sigmama).
\ee

%%%%%%%%%%%%%%%%%%%%%%%%%%%%%%%%%%%%%%%%%%%%%%%%%%%%%%%%%%%%%%%%%%%%%%%%%

%%%%%%%%%%%%%%%%%%%%%%%%%%%%%%%%%%%%%%%%%%%%

\section*{Acknowledgments}
Authors are grateful to Oleg Ogievetsky
and Rafael Stekolshchik for useful discussions.

%%%%%%%%%%%%%%%%%%%%%%%%%%%%%%%%%%%%%%%%%%%%%%%%%%%%%%%%%%%%%%%%%

\setcounter{equation}{0} \def\theequation{A\arabic{appen}.\arabic{equation}}

\newcounter{appen}
\newcommand{\appen}[1]{\par\refstepcounter{appen}
{\par\bigskip\noindent\large\bf Appendix \arabic{appen}. \medskip }{\bf \large{#1}}}

\renewcommand{\subsection}[1]{\refstepcounter{subsection}
{\bf A\arabic{appen}.\arabic{subsection}. }{\ \bf #1}}
\renewcommand\thesubsection{A\theappen.\arabic{subsection}}
\makeatletter \@addtoreset{subsection}{appen}

\renewcommand{\subsubsection}{\par\refstepcounter{subsubsection}
{\bf A\arabic{appen}.\arabic{subsection}.\arabic{subsubsection}. }}
\renewcommand\thesubsubsection{A\theappen.\arabic{subsection}.\arabic{subsubsection}}
\makeatletter \@addtoreset{subsubsection}{subsection}

%%%%%%%%%%%%%%%%%%%%%%%%%%%%%%%%%%%%%%%%%%%%%%%%%%%%%%%%%%%%%%%%%%%%

\appen{The proof of consistency condition (\ref{c1}) for $\lambda\neq \varkappa$}\label{appb}.

Let parameters $\mu_1  \stackrel{def}{=}  \mu^{K_1}$ and $\mu_2  \stackrel{def}{=}  \mu^{K_2}$ be such that
$\lambda_1  \stackrel{def}{=}  \lambda_{K_1}\neq \varkappa $ and $\lambda_2  \stackrel{def}{=}
\lambda_{K_2}\neq
\varkappa $. Let $b^1$ and $b^2$ denote $b_{K_1}$ and $b_{K_2}$ correspondingly.
Let us prove by induction that eqs. (\ref{c1}) are true.
To implement induction, we select a
part of degree $k$ in $\mu$
from eq. (\ref{nm1})  and observe that this part contains
a
degree $k$ polynomial in $b_M$ in the left-hand side
of eq. (\ref{nm1}) while the part on the right hand side of the
differential version (\ref{nm1}) of eq. (\ref{start}) which is of the same degree
in $\mu$ has degree $k-1$ as polynomial in $b_M$. This happens because of
the presence of the commutator $[b_L ,S]$ which is a zero degree polynomial
due to the basic relations (\ref{rel}). As a result, the inductive hypothesis allows
us to use the properties of $\varkappa$-trace provided that the above
commutator is always handled as the right hand side of eq. (\ref{rel}),
i.e., we are not
allowed to represent it again as a difference of the second-degree
polynomials.

Direct differentiation with the help of eq. (\ref{d}) gives
\bee
\label{p1}
(1-\varkappa \lambda_2 )\frac {\p}{\p\mu_2 }
\int
\,(-\varkappa \lambda_1 t_1 -t_2 ) \str \ig (
e^{t_1 S} [b^1,S]\,e^{t_2 S}\sigmama \ig ) D^1 t
- \ig (1 \leftrightarrow 2 \ig )
&=&\nn
=\left(\int \,
(1-\varkappa \lambda_2 )(-\varkappa \lambda_1 t_1 -t_2 )
\str \left( e^{t_1 S} [b^1 ,b^2 ]\,e^{t_2 S}
\sigmama  \right )D^1 t      \,
- \ig (1 \leftrightarrow 2 \ig )  \right)
&+&\nn
+\igg (
       \int (1-\varkappa \lambda_2 ) (-\varkappa \lambda_1 (t_1 +t_2 )-t_3 )
   \str \ig (
        e^{t_1 S} b^2 e^{t_2 S} [b^1 ,S]\,e^{t_3 S}
       \ig ) D^2 t\,\,
              - \ig (1 \leftrightarrow 2 \ig )
\igg )
&+&
\nn
+ \igg ( \int (1-\varkappa \lambda_2 )
          (-\varkappa \lambda_1 t_1 -t_2 -t_3 )
  \str  \ig(
          e^{t_1 S} [b^1 ,S]\,e^{t_2 S} b^2
          e^{t_3 S}\sigmama
      \ig )D^2 t \,
         - \ig (1 \leftrightarrow 2 \ig )
\igg)\,.&{}&
\eee

We have to show that the right hand side of eq. (\ref{p1}) vanishes.  Let us
first transform the second and the third terms on the right-hand side of eq.
(\ref{p1}). The idea is to move the operators $b^2$ through the exponentials
towards the commutator $[b^1 ,S]$ in order to use then the Jacobi identity for the
double commutators.  This can be done in two different ways inside  the
$\varkappa$-trace so that one has to fix appropriate weight factors for each of
these processes.  The correct weights turn out to be
\bee\label{p2}
D^2 t(-\varkappa \lambda_1 (t_1 +t_2 )-t_3 )b^2 \equiv
D^2 t(-\varkappa \lambda_1 -t_3 (1-\varkappa  \lambda_1 ))b^2=
\nn
D^2 t\left (\igg (\frac { \lambda_1 \lambda_2}{1-\varkappa \lambda_2}-t_3 (1-\varkappa
 \lambda_1 )\igg )
\overrightarrow {b^2} +
\frac {-\varkappa \lambda_1 }{1-\varkappa \lambda_2} \overleftarrow {b^2} \right )
\n
\eee
and
\bee\label{p3}
D^2 t(-\varkappa \lambda_1 t_1 -t_2 -t_3 )b^2
\equiv D^2 t ((-\varkappa \lambda_1 +1) t_1 -1)b^2=
\nn
D^2 t\left (
\igg (t_1 (1-\varkappa  \lambda_1 ) -\frac {1}{1-\varkappa \lambda_2} \igg )
\overleftarrow {b^2} -
\frac {-\varkappa \lambda_2 }{1-\varkappa \lambda_2} \overrightarrow {b^2} \right )
\n
\eee
in the second and third terms in the right hand side of eq. (\ref{p1}),
respectively.  Here the notation ${\overrightarrow A}$ and ${\overleftarrow
A}$ imply that the operator $A$ has to be moved from its position to the
right and to the left, respectively. Using eq. (\ref{r3})
along with the simple formula
\be\label{p4}
\int \, \phi (t_3 ,\ldots t_{n+1} )D^n t=
\int \, t_1\phi (t_2 ,\ldots t_n )D^{n-1} t
\n
\ee
we find that all terms which involve both $[b^1 ,S]$ and
$[b^2,S]$
cancel pairwise after antisymmetrization $1\leftrightarrow 2$.

As a result, one is left with some terms involving double commutators which
thanks to the Jacobi identities and antisymmetrization all reduce to
\be\label{1p5}
\int\,\ig (\lambda_1\lambda_2 t_1+t_2-t_1 t_2 (1-\varkappa \lambda_1)(1-\varkappa \lambda_2)\ig )
\str \ig (\exp (t_1 S)
[S,[b^1 ,b^2 ]]\exp (t_2 S) \sigmama \ig ) D^1 t\,.
\n
\ee
Finally, we observe that this expression can be equivalently rewritten in
the form
\be\label{p5}
\int \,
\ig (\lambda_1\lambda_2 t_1+t_2-t_1 t_2 (1-\varkappa \lambda_1)(1-\varkappa \lambda_2)\ig )
\left(\frac {\p}{\p t_1}-\frac {\p}{\p t_2} \right) \str \ig (\exp (t_1 S)
[b^1 ,b^2 ]\exp (t_2 S) \sigmama \ig )
D^1 t
\ee
and after integration by parts cancel the first term on the right-hand side
of eq. (\ref{p1}). Thus, it is shown that eqs. (\ref{nm1}) are mutually
compatible for the case $\lambda_{1,2}\neq \varkappa  $.

Analogously, we can show that eqs. (\ref{nm1}) are consistent with eq.
(\ref{m1}).  Indeed, let $\lambda_1 =\varkappa $, $\lambda_2 \neq \varkappa $. Let us prove
that
\be\label{p6}
\frac {\p}{\p\mu_2}\str \ig ([b^1 , \exp (S)]\sigmama  \ig )=0
\ee
provided that the $\varkappa$-trace is well-defined for the lower degree polynomials.
The
explicit differentiation gives
\bee\label{p7}
\!\!\!\!\!\!\!\!\!\!\!\!\!\!\!\!\!\!
\frac {\p}{\p\mu_2}\str \ig ([b^1 , \exp (S)]\sigmama  \ig )&=&
\int \,\str \ig ([b^1 , \exp (t_1 S)b^2  \exp (t_2 S)]\sigmama  \ig
)D^1 t
\nn &=&
(1-\varkappa \lambda_2 )^{-1} \str \ig ([b^1 ,(b^2
 \exp (S) -\varkappa  \lambda_2  \exp ( S)b^2 )]\sigmama  \ig ) +\ldots
\eee
where dots denote some terms of the form $\str \ig( [b^1 , B]\sigmama \ig )$
involving more commutators inside $B$, which therefore amount to some
lower degree polynomials and vanish by the inductive hypothesis. As a result,
we find that
\bee\label{p8}
\!\!\!\!\!\!\!\!\!\!\!\!\!\!\!\!\!\!
\frac {\p}{\p\mu_2}\str \ig ([b^1 , \exp (S)]\sigmama  \ig )&=&
(1-\varkappa \lambda_2 )^{-1} \str \ig ((b^2 [b^1 , \exp (S)] -\varkappa  \lambda_2
[b^1 , \exp ( S)]b^2 )\sigmama  \ig )\nn
&+& (1-\varkappa \lambda_2 )^{-1} \str \ig (([b^1 ,b^2 ] \exp (S) -\varkappa
\lambda_2  \exp ( S)[b^1 ,b^2 ])\sigmama  \ig )\,.
\eee
This expression vanishes by the inductive hypothesis, too.

\newpage
\appen{The proof of consistency conditions (\ref{c3}) (the case of special polynomials)}\label{appc}

In order to prove eq. (\ref{c3}) we use the inductive hypothesis {\bf (i)}.
In this appendix we use  the convention that any
expression with the coinciding
upper or lower indices
are automatically symmetrized, e.g.,
$F^{II} \stackrel{def}{=} \frac  1 2 (F^{I_1 I_2}+F^{I_2 I_1})$.
Let
us write the identity
\be\label{p12}
0 = \sum_{M}\str \ig (\ig [  \exp (S^\prime )
\{ b_I ,f_{IM}\}\mu^M ,b_J b_J\ig ]
\sigmama \ig ) -(I \leftrightarrow J)
\ee
which holds
due to Lemma \ref{lemma4}
for all terms of degree $k-1$ in $\mu$
with $E(g) \leqslant l+1$
and
for all lower degree polynomials in $\mu$
(one
can always move $f_{IJ}$ to $\sigmama $ in eq. (\ref{p12}) combining $f_{IJ}\sigmama$
into a combination of
elements of $W(\R)$ analyzed in Lemma \ref{lemma4}).

Straightforward calculation of the commutator in the right-hand-side of eq.
(\ref{p12}) gives $0 = X_1+X_2+X_3$, where
\bee\label{1p13}
X_1&=&-\sum_{M,L}\int\,\str
          \left (   \exp (t_1 S^\prime ) \{b_J , F_{JL} \}
              \mu^L  \exp (t_2 S^\prime ) \{b_I ,f_{IM}\}\mu^M \sigmama
          \right) D^1 t
       -(I \leftrightarrow J)\,,\nn
X_2&=&\sum_M \str
          \left (  \exp (S^\prime )
              \ig \{\{b_J, F_{IJ}\},f_{IM}
              \ig\}  \mu^M \sigmama
          \right)-(I \leftrightarrow J)\,,\nn
X_3&=& \sum_M \str
              \left (  \exp (S^\prime )
                     \ig
\{b_I,\{b_J,[f_{IM},b_J]\}
                     \ig \}\mu^M \sigmama
               \right)-(I \leftrightarrow J)\,.
\eee
The terms bilinear in $f$ in $X_1$ cancel due to the
antisymmetrization ($I\leftrightarrow J$) and the inductive hypothesis
{\bf (i)}. As a result, one
can transform $X_1$ to the form
\bee\label{p13}
X_1
=
\left (- \frac  1 2
\left [ L_{JJ},\,R_{II}\right ] +2
\str \ig ( e^{S^\prime }
\{b_I ,f_{IJ}\}\mu_J \sigmama \ig )
\right)
-
(I \leftrightarrow J).
\eee

Substituting
$F_{IJ}={\cal C}_{IJ}+ f_{IJ}$ and $f_{IM}= ([b_I,b_M] -{\cal C}_{IM})$
one transforms $X_2$ to the form
\bee\label{x2}
X_2
&=& 2{\cal C}_{IJ}R_{IJ}
-2 \left( \str
         \ig ( e^{S^\prime } \{b_J ,f_{IJ}\}\mu_I \sigmama
         \ig ) -(I \leftrightarrow J)
   \right)+Y,
\eee
where
\bee\label{Y}
Y=  \str \left ( e^{S^\prime }
              \ig \{ \{b_J, f_{IJ}\},[b_I,\,S^\prime]
              \ig\}   \sigmama
          \right)
    -(I \leftrightarrow J)\,.
\eee
Using that
\bee\label{S}
 \str
    \left (  \exp (S^\prime )
       \left[   P f_{IJ}Q,\,S^\prime
       \right]   \sigmama
    \right)
=0
\eee
provided that the inductive hypothesis can be used,
one transforms $Y$ to the form
\bee\label{x2x2}
Y=
\str
 \bigg ( e^{S^\prime }
   \bigg( &-&[f_{IJ},\, (b_I S^\prime b_J + b_J S^\prime b_I )]
             {} -
                b_I [f_{IJ},\, S^\prime] b_J -b_J [f_{IJ},\, S^\prime] b_I\nn
              & + & [f_{IJ},\, \{b_I\,, b_J\}] S^\prime
    \bigg) \sigmama
 \bigg).
\eee

Let us rewrite $X_3$ in the form $X_3=X_3^{s}+X_3^{a}$, where
\bee
X_3^s=
\frac  1 2 \sum_M \str
              \left ( e^{S^\prime }
               \ig (
                    \ig\{b_I,\{b_J,[f_{IM},b_J]\}\ig \} +
                    \ig\{b_J,\{b_I,[f_{IM},b_J]\}\ig \}
                \ig )
                          \mu^M \sigmama
               \right)
-
(I \leftrightarrow J)\,,\nn
X_3^a=
\frac  1 2 \sum_M \str
              \left ( e^{S^\prime }
               \ig (
                    \ig\{b_I,\{b_J,[f_{IM},b_J]\}\ig \} -
                    \ig\{b_J,\{b_I,[f_{IM},b_J]\}\ig \}
                \ig )
                          \mu^M \sigmama
               \right)
-
(I \leftrightarrow J)\,. \nonumber
\eee
With the help of the Jacobi identity
$[f_{IM},b_J]-[f_{JM},b_I]=[f_{IJ},b_M]$
one expresses $X_3^s$ in the form
$$
X_3^s=
\frac  1 2 \str
\left ( e^{S^\prime}
 \left (
    \{b_I,b_J\}[f_{IJ}, S^\prime] + [f_{IJ}, S^\prime]\{b_I,b_J\}
     +2 b_I[f_{IJ}, S^\prime]b_J + 2 b_J[f_{IJ}, S^\prime]b_I
  \right ) \sigmama
\right ).
$$
Let us transform this expression for $X_3^a$ to the form
\bee\label{100}
X_3^a=
\frac  1 2 \sum_M
\str
\left ( e^{S^\prime}
      \left[  F_{IJ},\, [f_{IM},b_J]\right ] \mu^M \sigmama
\right )       -(I \leftrightarrow J).
\eee

Substitute
$F_{IJ} ={\cal C}_{IJ}+ f_{IJ}$ and $f_{IM}= ( [b_I,b_M]-{\cal C}_{IM})$
in eq. (\ref{100}).
After simple transformations we find that $Y+X_3=0$.
From eqs. (\ref{p13}) and (\ref{x2}) it follows that
the right hand side of eq. (\ref{p12}) is
equal to
$$
\frac  1 2 ([L_{II},\, R_{JJ}]-[L_{JJ},\, R_{II}]) + 2{\cal C}_{IJ} R_{IJ}.
$$
This completes the proof of the consistency conditions
(\ref{c3}).

%%%%%%%%%%%%%%%%%%%%%%%%%%%%%%%%%%%%%%%%%%%%%%%%%%%%%%%%%%%%%%%%%

\end{document}